\RequirePackage{fix-cm}
\documentclass[11pt]{article} 
\usepackage[margin=2.5cm]{geometry}

\usepackage{multirow}
\usepackage{mathtools}
\usepackage{subfigure}
\usepackage{natbib}
\usepackage{booktabs}

\DeclareMathOperator{\ddiv}{div}

\newcommand{\ed}{\ensuremath{\mathrm{d}}} 
\newcommand{\mc}[1]{\ensuremath{\mathcal{#1}}}

\newcommand{\mr}[1]{\ensuremath{\mathrm{#1}}}
\newcommand{\der}[2]{\ensuremath{\frac{\ed #1}{\ed #2}}}

\usepackage{authblk}

\title{Scale-selective dissipation in energy-conserving finite element schemes for two-dimensional turbulence}
\author{Andrea Natale and Colin J. Cotter}
\affil{Department of Mathematics, Imperial College London, London, SW7 2AZ, UK}
\date{\today}
\newcommand{\col}[1]{#1} %{{\color{blue}#1}}
\usepackage{bm}
\newcommand{\mm}[1]{{\bm{#1} }}

\begin{document}

%\runningheads{A. Natale and C. J. Cotter}{Scale-selective dissipation in energy-conserving finite elements}

\maketitle
\begin{abstract}
We analyse the multiscale properties of energy-conserving upwind-stabilised finite element discretisations of the two-dimensional incompressible Euler equations.
We focus our attention on two particular methods: the Lie derivative discretisation introduced in \cite{natale16} and the \col{Streamline Upwind/Petrov-Galerkin (SUPG)} discretisation of the vorticity advection equation.
Such discretisations provide control on enstrophy by modelling different types of scale interactions.
We quantify the performance of the schemes in reproducing the non-local energy backscatter that characterises two-dimensional turbulent flows.
\end{abstract}

\section{Introduction}
% Motivation behind the work, physics of 2D turbulent flow.
% Numerical methods and performance study of Thuburn et al
% A new upwind form that conserves energy.
% We show that this gives rise to an implicitely multiscale.
% We study the performace as in Thuburn

Numerical simulations are usually unable to fully represent all the dynamically relevant scales in atmospheric flows \citep{williamson07,shutts05,thuburn14}. 
%In modern dynamical cores (mention grid size and relevant scales size). [...] 
For such unresolved flows, different discretisations should be compared not only in terms of accuracy and order of convergence, but also on the way the unresolved scales affect the computation. 
%This is a particularly pressing issue given the scale separation that characterises atmospheric flows.
More specifically, as both kinetic energy and enstrophy are statistically relevant in determining large scale flow evolution \citep{majda06}, their dynamics within the resolved scales should be appropriately captured by the numerical scheme in use. 

At low resolution, energy conservation is generally considered an appropriate requirement for numerical schemes modelling atmospheric flows \citep[see, e.g.,][]{shutts05}.
Enstrophy conservation, on the other hand, induces a piling up of enstrophy at small scales, a phenomenon known as ``spectral blocking". 
Such a phenomenon is well known, and it has been studied extensively in the context of two-dimensional turbulence.
The classical picture is described in \cite{maltrud93}, where the different behaviour of energy and enstrophy transfer are studied in terms of triads of wavenumbers.
In particular, one sees that in the enstrophy intertial range the most significant triads interactions are very nonlocal, with two of the three wavenumbers being much higher than the third. 
The large wavenumbers are dominant when considering enstrophy transfer, which can then be considered to be local.
Energy transfer, on the other hand, also has a non negligible nonlocal contribution due to the smallest wavenumber, which causes an upscale transfer known as backscatter.

The general approach in numerical schemes is to focus on enstrophy dissipation, to avoid its piling up at small scales.
This however ignores the energy backscatter and yields simulations that are too dissipative.
One possible solution is achieved by combining some form of artificial hyper-viscous dissipation with an ``energy fixer'' \citep{shutts05}, which reinserts the dissipated energy in the form of random perturbations or specific large-scale velocity patterns. 

\cite{thuburn14} proposed an energy fixer modelled on the instantaneous vorticity distribution and applied it to various finite difference and finite volume schemes for the barotropic vorticity equation, i.e.\ the two-dimensional incompressible Euler equations, which they  considered as a model problem. 
Moreover, they compared the behaviour of the various schemes in terms of energy and enstrophy tendencies in spectral space. 
In particular, they quantified how well the discretisations would reproduce the energy and enstrophy tendencies due to small scales in a truncated DNS simulation.

In this paper, we perform a similar analysis, but apply it specifically to finite element upwind-stabilised and energy-conserving discretisations.
The reason for our interest is that upwind stabilisation in finite element methods can be regarded as a form of scale-selective dissipation, providing a way of modelling the interaction of scales in the simulation. 
Depending on the scheme, such an interaction can be either between resolved and quasi-resolved scales, or between resolved and unresolved scales (as a proper subgrid model).
This property has been used extensively in designing discretisations for fluid simulation \citep[see, e.g.,][]{hughes95,burman07,becker11}, but to the authors' knowledge it has never been analysed in the context of energy conserving simulations. 
In recent years, finite element methods have been developed and considered for
numerical weather prediction and climate modelling
\citep{fournier2004spectral,thomas2005ncar,dennis2011cam,kelly2012continuous,kelly2013implicit,brdar2013comparison,bao2015horizontally,marras2015review}.
The use of stabilised finite element methods in atmospheric simulations
has been investigated in \citep{marras2013simulations}. Compatible finite
element methods have been developed in order to extend conservation and stability properties from C-grid staggered finite difference methods to the finite element setting \citep{cotter2012mixed,mcrae2014energy,natale2016compatible}; these methods provide the context for this paper.

Here, we focus our attention on two algorithms that model different scales interactions while preserving kinetic energy. 
The first is
the discretisation developed in \cite{natale16}, which consists of a finite element ${H}(\ddiv)$-conforming scheme with an energy-conserving upwind stabilisation. Such a discretisation was derived using Hamilton's principle of stationary action applied to an appropriately defined Lagrangian, representing the total kinetic energy. Then, the scheme conserves the Lagrangian by construction and possesses the attractive feature of a variational derivation, with possible extensions to different Hamiltonian models. Furthermore, we will observe that the conservative upwind stabilisation implicitly induces an energy transfer from small quasi-resolved scales to fully resolved scales, avoiding in this way spectral blocking without dissipating energy.
Our second focus will be 
the SUPG stabilised discretisation of the vorticity/stream function formulation of the Euler equations introduced in \cite{tezduyar88,tezduyar89}. The SUPG stabilisation can be interpreted as a subgrid model in which the unresolved scales are proportional to the residual at each time step \citep{hughes95}. Moreover, as the stabilisation is applied to the vorticity advection equation, only enstrophy is dissipated while kinetic energy is conserved exactly.

This paper is organised as follows. In Section~\ref{sec:fem} we introduce the finite element framework and we describe the different schemes analysed in the paper. In Section~\ref{sec:multi} we discuss the multiscale interpretation of the upwind stabilisation. In Section~\ref{sec:num} we present some numerical tests demonstrating the behaviour of the schemes in terms of spectral energy and enstrophy tendencies. Conclusions are presented in Section~\ref{sec:con}.

% We only use Euler.
% Stream function issue

\section{Finite element framework}\label{sec:fem}
In this section we introduce some basic finite elements tools (Section~\ref{sec:spaces}) and describe the methods studied in this paper (Section~\ref{sec:methods}). In particular, as we treat mixed formulations in velocity/pressure as well as stream function/vorticity form, we will need to introduce finite element spaces for all these variables.
%Note that in two dimensions one could express the mixed velocity/pressure scheme of this section, in a completely equivalent way, in terms of stream function only. This is because we use ${H}(\ddiv)$-conforming spaces for the velocity which imply exact enforcement of the divergence-free condition \citep{Boffi13}. 
%However, formulating the problem in terms of the velocity field gives more insight to the multiscale character of the scheme. 

\subsection{Preliminaries}\label{sec:spaces}
Consider a two-dimensional domain $\Omega$. We denote by $\mc T$ a triangulation of $\Omega$, i.e.\ a decomposition of $\Omega$ in triangular elements $T\subset \Omega$. %We define the Continuous Galerkin (CG) finite element space as follows,
%\begin{equation}
%CG_r({\mc T})\coloneqq \{ \psi\in C^0(\Omega) \,:\, \psi|_T \in {\mc P}_r(T)  \}\,.
%\end{equation}
%In other words, $CG_r({\mc T})$ is the space of continuous functions that are polynomials of degree up to $r$ in each element. Furthermore, we define
%\begin{equation}
%\overline{CG}_r({\mc T})\coloneqq \left\{ \psi\in CG_r(\mc T) \,:\, \int_\Omega \psi \,\ed x =0  \right\}\,.
%\end{equation}
%In the following, we will define the stream function to be an element of $V\coloneqq\overline{CG}_r(\mc T)$ for a fixed $r$, with $\Omega$ being a periodic square domain. As a consequence, the velocity field $u \coloneqq \nabla^\perp \psi$ is not necessarily a continuous vector field. 
%The velocity discontinuities are taken into account in the discretisation via the jump and average operators, which we now define. 
Let $\mc E$ be the set of element edges of the triangulation $\mc T$. We fix an orientation for all edges $e\in {\mc E}$ by specifying the unit normal vector ${\mm n}_e$. An orientation on the edge $e\in {\mc E}$ defines a positive and negative side on $e$, so that any finite element scalar or vector function $f$ on $\Omega$ has two possible restrictions on $e$, denoted $f^+$ and $f^-$, with $f^+\coloneqq f|_{T^+}$ and $T^+\in \mc T$ being the element with outward normal ${\mm n}_e$. Therefore, we can define on each $e\in \mc E$ the jump operator  $[f] \coloneqq f^+ - f^-$ and the average operator $\{f\} \coloneqq (f^+ + f^-)/2$. 

We start with defining the stream function and vorticity space. For any element $T\in \mc T$, we denote by ${\mc P}_r(T)$ the linear space of polynomials on $T$ of degree up to $r$. Then the stream function and vorticity space is given by
\begin{equation}
\begin{aligned}
V_0 = CG_{r+1} \coloneqq \{ \psi \,:\, \psi|_T \in & {\mc P}_{r+1}(T) ~~\forall\, T\in {\mc T},\\& \psi^+|_e = \psi^-|_e ~~\forall\, e\in {\mc E} \},
\end{aligned}
\end{equation}
which is the Continuous Galerkin finite element space of degree $r+1$. Clearly, functions in $V_0$ are continuous on $\Omega$, i.e.\ $V_0\subset C^0(\Omega)$. More precisely, as we will employ periodic boundary conditions, the stream function space will be the space $\bar{V}_0$ of functions in $V_0$ with zero average on $\Omega$.

We now define the velocity and pressure spaces.
These are defined in order to produce a stable mixed finite element discretisation, see \cite{Boffi13} for details.  The velocity space ${\mm V}_1$ is given by
\begin{equation}
\begin{aligned}
{\mm{V}}_1 = {\mm{BDM}}_r \coloneqq \{ {\mm u} \,:\, & {\mm u}|_T \in ({\mc P}_r(T))^2 ~~\forall\, T\in {\mc T},\\ & {\mm u}^+|_e \cdot {\mm n}_e = {\mm u}^-|_e \cdot {\mm n}_e ~~\forall\, e\in {\mc E} \},
\end{aligned}
\end{equation}
which is the Brezzi-Douglas-Marini finite element space of degree $r$ \citep{Boffi13}. Note that a vector field in ${\mm V}_1$ is not completely continuous, as only the normal components at the edges are continuous. 
The subspace of ${\mm V}_1$ of continuous velocity fields is denoted by ${\mm V}_1^l$, i.e.\ ${\mm V}_1^l\coloneqq {\mm V}_1 \cap (C^0(\Omega))^2$, and it will be used to model the large scale component of the velocity. The quasi-resolved small scale velocity space ${\mm V}_1^s$ is the orthogonal complement of ${\mm V}_1^l$ in ${\mm V}_1$, with respect to the ${\mm L}^2$ inner product. In other words, we have ${\mm V}_1 = {\mm V}_1^l \oplus {\mm V}_1^s$.

The pressure space $V_2$ is given by
\begin{equation}
V_2 = DG_{r-1} \coloneqq \{ p \,:\, p|_T \in {\mc P}_{r-1}(T) ~~\forall\, T\in {\mc T}\},
\end{equation}
which is the Discontinuous Galerkin finite element space of degree $r-1$ \citep{Boffi13}. More precisely, the pressure space will be the space $\bar{V}_2$ of functions in $V_2$ with zero average on $\Omega$.

For any $r$, the couple $({\mm V}_1,\bar{V}_2)$ can be used to produce a stable mixed finite element formulation and enforce the divergence-free constraint exactly on the velocity field. In practice, this constraints the velocity to the space $\nabla^\perp V_0 \subseteq {\mm V}_1$, where the superscript $\perp$ denotes a clockwise rotation of $\pi/2$.
In principle, other choices for the spaces ${\mm V}_1$ and $\bar{V}_2$ are possible, however, for simplicity, we will restrict ourself to the spaces mentioned in this section. 

Finally, for any scalar fields $f$ and $g$ on $\Omega$, we define $(f,g)_\Omega \coloneqq \int_\Omega f \cdot g \, \ed x$ and $\|f \|_\Omega \coloneqq \sqrt{(f,f)_\Omega}$. The same notation is used for vector and tensor fields, and for integrals over edges or elements.

%Finite element
%Domain and decomposition
%Continuous Galerkin finite element space (stream function space)
%Velocity field
%Jumps and avarages
%Approximation properties

\subsection{Upwind methods for the two-dimensional Euler equations}\label{sec:methods}

We now introduce two categories of upwind-stabilised finite element schemes for the two-dimensional Euler equations. The first is based on the velocity/pressure formulation and the second on the stream function/vorticity formulation. We start with the former. The Euler equations can be written in terms of velocity and pressure as follows, 
\begin{equation} \label{eq:euler}
\left\{
\begin{array}{l}
\displaystyle \partial_t{\mm u} + ({\mm u} \cdot \nabla) {\mm u} + \nabla p = 0, \\
\nabla\cdot {\mm u} =0,
\end{array}
\right.
\end{equation}
where ${\mm u}$ and $p$ are the velocity and pressure fields respectively, and the domain $\Omega$ is assumed to be doubly-periodic. Using the divergence-free condition the first equation in the system \eqref{eq:euler} can be rewritten as
\begin{equation}\label{eq:fluxform}
\displaystyle \partial_t{\mm u} + \nabla \cdot({\mm u} \otimes {\mm u}) + \nabla p = 0\, ,
\end{equation} 
which we refer to as flux form. Alternatively, using standard vector calculus identities, we can also write it as
\begin{equation}\label{eq:lieform}
\displaystyle \partial_t{\mm u} + {\mm u}^\perp (\nabla^\perp\cdot {\mm u}) + \nabla P = 0\, ,
\end{equation}
where $P\coloneqq p + \|{\mm u}\|_2^2/2$, \col{and $\|\cdot\|_2$ denotes the Euclidean norm}. We refer to Equation~\eqref{eq:lieform} as Lie derivative formulation, as the advection term in this formulation may be interpreted as the Lie derivative of the velocity one-form \citep{natale16}.
Then, the general form of the mixed finite element discretisation of the system in \eqref{eq:euler} requires us to find $({\mm u},\tilde{p})\in {\mm V}_1\times \bar{V}_2$ such that
\begin{equation} \label{eq:eulerw}
\left\{
\begin{array}{ll}
\displaystyle (\partial_t{\mm u},{\mm v} )_\Omega + a({\mm u};{\mm u},{\mm v}) + s({\mm u};{\mm u},{\mm v}) -(\tilde{p}, \nabla\cdot {\mm v})_\Omega = 0, \\
(\nabla\cdot {\mm u}, q)_\Omega =0,
\end{array}
\right.
\end{equation}
for all $({\mm v},q)\in {\mm V}_1 \times \bar{V}_2$, where $\tilde{p}$ can represent either $P$ or $p$. For our choice of finite element spaces, the system \eqref{eq:eulerw} ensures that $\nabla \cdot {\mm u} = 0$ is satisfied pointwise at all times. The forms $a(\cdot; \cdot,\cdot)$ and $s(\cdot;\cdot,\cdot)$ define the discretisation of the advection term, and in particular the second defines the upwind stabilisation. 
We consider here two particular choices.

\paragraph{Flux form discretisation} This is the standard discretisation used in DG schemes for the Navier-Stokes equations \citep[see, e.g.,][for an application to the Euler equations]{Cockburn04,Guzman15}.
Note that such discretisation is not energy-conserving, and therefore we will only use it as a reference due to its similar structure to the Lie derivative discretisation, introduced below. The flux form scheme is derived by writing the advection term in flux form as in Equation~\eqref{eq:fluxform} and integrating by parts.
It is defined by
\begin{align}\label{eq:dg_a}
\begin{split}
a({\mm u};{\mm u},{\mm v}) \coloneqq & -\sum_{T\in {\mc T}}({\mm u} , ({\mm u} \cdot \nabla)  {\mm v} )_T\\& + \sum_{e\in {\mc E}} ({\mm u}\cdot {\mm n}_e\,\{{\mm u}\},[{\mm v}])_e \, ,
\end{split}\\ \label{eq:dg_s} s({\mm u};{\mm u},{\mm v}) \coloneqq & \sum_{e\in {\mc E}} ( c_e\, {\mm u} \cdot {\mm n}_e \, [{\mm u}],[{\mm v}])_e\,,
\end{align}
where $c_e$ is a function defined on the edges, and dependent on $u$. Generally, $c_e = \alpha\, {\mm u}\cdot {\mm n}_e/(2|{\mm u}\cdot {\mm n}_e|)$, where $\alpha>0$ is a constant defining the level of upwinding. For $\alpha=1$, the two terms can be combined by replacing the average in the advection form $a$ with the upwind value of ${\mm u}$, i.e.\ ${\mm u}^+$ if ${\mm u}\cdot {\mm n}_e>0$ and ${\mm u}^-$ otherwise. For $0<\alpha<1$, the average term becomes a skewed average. 

Note that, if $c_e = \alpha\, {\mm u}\cdot {\mm n}_e/(2|{\mm u}\cdot {\mm n}_e|)$ then $s({\mm u};{\mm u},{\mm u})\geq 0$, therefore it represents a dissipation term in the discretisation.  Denoting by $h$  the maximum element diameter in the triangulation $\mc T$, we define $\tau_h\coloneqq h/{\sup_\Omega (\|{\mm u}\|_2)}$ to be the local time scale associated to the mesh. Then, the dissipation induced by the upwinding acts on a time scale bounded from above by $\tau_h/\alpha$. This follows immediately from the following string of inequalities,
\begin{equation}
\begin{aligned}
  s({\mm u};{\mm u},{\mm u}) &\leq  \frac{\alpha}{2} \sup_\Omega (\|{\mm u}\|_2) \sum_{e\in {\mc E}}  \, ([{\mm u}],[{\mm u}])_e, \\
  &\leq C \alpha \sup_\Omega (\|{\mm u}\|_2)h^{-1} \|{\mm u}\|^2_\Omega =  C \alpha\tau_h^{-1} \|{\mm u}\|^2_\Omega,
\end{aligned}
\end{equation}
where $C>0$ is a constant independent of $h$.
\paragraph{Lie derivative discretisation} We now consider the discretisation introduced in \cite{natale16}. It can be derived by writing
the advection term in Lie derivative form as in Equation~\eqref{eq:lieform} and integrating by parts. It is defined by
\begin{align}\label{eq:conservative_a}
\begin{split}
a({\mm u};{\mm u},{\mm v}) \coloneqq  & \sum_{T\in {\mc T}}({\mm u}^\perp, \nabla ( {\mm u}^\perp \cdot {\mm v}) )_T \\&- \sum_{e\in {\mc E}} (\{{\mm u}^\perp\}\cdot n_e,[{\mm u}^\perp \cdot {\mm v}])_e \, ,
\end{split}\\
\label{eq:conservative_s}
s({\mm u};{\mm u},{\mm v}) \coloneqq & -\sum_{e\in {\mc E}} (c_e [{\mm u}^\perp]\cdot n_e,[{\mm u}^\perp \cdot {\mm v}])_e\, .
\end{align}
Note that in this case $s({\mm u};{\mm u},{\mm u})=0$ independently of the choice of the function $c_e$. Such a property directly leads to conservation of energy. This can be verified by setting ${\mm v}={\mm u}$ and $q=\tilde{p}$ in Equation~\eqref{eq:eulerw}, and noting that we also have $a({\mm u};{\mm u},{\mm u})=0$. However, interpreting the stabilisation term as we did for the flux form discretisation is not straightforward. \cite{natale16} noticed that for smooth advecting velocity such a discretisation coincides with the Eulerian Lie derivative discretisation proposed in \cite{Heumann11} for the linear advection-diffusion problem. In Section~\ref{sec:multi} we elaborate on this observation, and we give a multiscale interpretation to the stabilisation term.
\\

%Include convergence estimates?
Numerical tests show that both the flux form and the Lie derivative scheme converge with order $r+1$, in terms of velocity ${\mm L}^2$ error, when standard upwinding is employed, i.e.\ $\alpha =1$, see \cite{Guzman15} and \cite{natale16}. A priori convergence estimates with suboptimal convergence rate $r$ can be found in the same references.

We now consider the stream function/vorticity formulation of the two-dimensional Euler equations which is the basis for the SUPG discretisation. This is given by
\begin{equation} \label{eq:eulervort}
\left\{
\begin{array}{l}
\displaystyle \partial_t{\omega} + {\mm u} \cdot \nabla \omega  = 0,\\
\Delta \psi = \omega,
\end{array}
\right.
\end{equation}
where ${\mm u}\coloneqq \nabla^\perp \psi$. The general form of the mixed finite element discretisation of the system in \eqref{eq:euler} is given by : find $(\psi,\omega)\in \bar{V}_0\times{V}_0$ such that
\begin{equation} \label{eq:eulervw}
\left\{
\begin{array}{ll}
\displaystyle (\partial_t{\omega},\varphi )_\Omega + \tilde{a}(u;\omega,\varphi) + \tilde{s}(u;\omega,\varphi)  = 0, \\
(\nabla\psi, \nabla \phi)_\Omega = -(\omega,  \phi)_\Omega,
\end{array}
\right.
\end{equation}
for all $(\varphi,\phi) \in V_0 \times \bar{V}_0$.
\paragraph{SUPG discretisation} The SUPG discretisation is based on the standard Galerkin discretisation for vorticity advection, supplemented with a stabilisation term  which is taken to be proportional to the residual of the vorticity advection equation  \citep{tezduyar88,tezduyar89}, i.e.\ we choose
\begin{align}\label{eq:supg_a}
\tilde{a}({\mm u};\omega,\varphi) &\coloneqq  ({\mm u}\cdot \nabla \omega, \varphi)_\Omega, \, \\
\label{eq:supg_s}
\tilde{s}({\mm u};\omega,\varphi) &\coloneqq   (\tau_{s} R(\omega), {\mm u}\cdot \nabla \varphi)_\Omega,
\end{align}
where $R(\omega) \coloneqq \partial_t \omega+{\mm u} \cdot \nabla \omega$ and where $\tau_s\geq0$ is a scalar function that may be discontinuous across elements.  Remarkably, as noted in \citet{mcrae15}, the SUPG scheme still preserves energy regardless of the stabilisation, just as for the Lie derivative scheme. This can be verified directly by noting that $\tilde{a}({\mm u};\omega,\psi) = \tilde{s}({\mm u};\omega,\psi)=0$ and setting $\varphi = \psi$ in Equation~\eqref{eq:eulervw}. The stabilisation term is trivially consistent being proportional to the residual. Moreover it induces streamline dissipation on the enstrophy density. %Note that, alternatively, we can obtain this discretisation by modifying the test function in the standard Galerkin discretisation, i.e. taking as test function $\tilde{\varphi} = \varphi + \tau_{s} u \cdot \nabla \varphi$.
The factor $\tau_s$ defines the time scale at which the dissipation acts and it is referred to as ``intrinsic time scale" \citep{hughes95}. Generally speaking, it can be taken to be proportional to the global value $\tau_h$. However, defining $\tau_s=\tau_h$ may lead to over-dissipative solutions, especially for higher order elements \citep{almeida97,codina92}. For the linear advection problem with high order elements, \cite{almeida97} suggest the value
\begin{equation}\label{eq:tau}
\tau_s = \frac{  \beta h_T \xi}{2 \|{\mm u}\|_2 }
\end{equation}
where $h_T$ is the characteristic length of the element $T\in \mc T$, depending on the shape and the local velocity direction of the element, $\xi=1$ for linear elements and $\xi=1/2$ for quadratic elements.  The coefficient $\beta$ is a constant; in \cite{almeida97} $\beta=1$, but we will also consider different values in our numerical tests.

We note that the stream function/vorticity SUPG scheme can be
re-interpreted as a mixed velocity-pressure scheme, since if $\psi\in
V_0$, then $\nabla^\perp\psi$ spans the divergence-free subspace of ${\mm V}_1$,
and in the absence of boundaries, $\omega\in V_0$ can be obtained by solving
\begin{equation}
  (\gamma,\omega)_\Omega = -(\nabla^\perp\gamma, {\mm u})_\Omega, 
\end{equation}
for all test functions $\gamma$ in $V_0$.
Then we obtain an equivalent formulation that is an approximation
of the incompressible Euler equations in Lie derivative form, with
\begin{align}
a({\mm u};{\mm u},{\mm v}) =  & ({\mm v},{\mm u}^\perp\omega)_\Omega, \\
s({\mm u};{\mm u},{\mm v}) = & ({\mm v},-{\mm u}^\perp\tau_s R(\omega))_\Omega.
\end{align}
This formulation allows us to extend the energy-conserving SUPG
approach to compatible finite element method discretisations of the
shallow water equations in velocity-height formulation, following
\citet{mcrae2014energy}. \col{The equivalent discretisation is
\begin{align}
  (\partial_t\mm{u},\mm{v})_{\Omega} + ({\mm v},{\mm F}^\perp q)_\Omega \qquad \qquad & \nonumber \\
 \nonumber + ({\mm v},-{\mm F}^\perp \tau_s \left(q_t+\mm{u}\cdot \nabla q\right))_\Omega & \\  - \left(\nabla\cdot\mm{v},gh + \frac{1}{2}|\mm{u}|^2\right)_\Omega & = 0, \\
  (\partial_t h + \nabla \cdot \mm{F},\phi)_\Omega & = 0, \\
  (\mm{F},\mm{w})_\Omega - (\mm{u}h,\mm{w})_\Omega & = 0, \\
  (\gamma,qh)_\Omega + (\nabla^\perp\gamma, u)_\Omega - (\gamma, f)_\Omega & = 0, 
\end{align}
for all test functions $(\mm{v},\phi,\mm{w},\gamma)\in {\mm V}_1\times V_2\times {\mm V}_1\times V_0$, having introduced potential vorticity $q\in V_0$ and mass flux
$\mm{F}\in {\mm V}_1$, and where $g$ is the acceleration due to gravity and
$f$ is the Coriolis parameter.}

%To the authors' knowledge, a priori estimates for the SUPG discretisations are not available in literature. However, the unstabilised Galerkin discretisation (i.e. the case $\tau_s=0$) converges with order $r-1/2$ in $L^2$ for both the vorticity and for the velocity field \citep[see][]{Liu01}.

%Definition of the upwind factor for quadratic simplicial elements \cite{codina92} or higher order global upwind coefficient in the element \cite{almeida97}.  

%Interpretation as sub-grid model \cite{codina11}
%Origin of SUPG described in for compressible Euler \cite{hughes10} .

\section{Multiscale interpretation}\label{sec:multi}
We now describe the multiscale character of the Lie derivative (Section~\ref{sec:liescale}) and the SUPG discretisations (Section~\ref{sec:supgscale}). We show that the Lie derivative scheme models the interaction of resolved and quasi-resolved scales within the velocity space ${\mm V}_1$, whereas the SUPG scheme models the interaction of the resolved and unresolved scales in the vorticity advection equation.

\subsection{Resolved/quasi-resolved interactions in the Lie derivative discretisation}\label{sec:liescale}
Consider the Lie derivative scheme and recall the decomposition ${\mm V}_1 = {\mm V}_1^l\oplus {\mm V}_1^s$ defined in Section~\ref{sec:spaces} . Let ${\mm u}={\mm u}^l+{\mm u}^s$, where ${\mm u}^l$ and ${\mm u}^s$ are the ${\mm L}^2$ projections of ${\mm u}$ onto ${\mm V}_1^l$ and ${\mm V}_1^s$ respectively. Then, taking ${\mm v}={\mm u}^l$ in the expression for $s$ in Equation~\eqref{eq:conservative_s},  we obtain
\begin{equation}
\begin{aligned}
s({\mm u};{\mm u},{\mm u}^l) &=  -\sum_{e\in {\mc E}} (c_e [{\mm u}^\perp]\cdot {\mm n}_e,[{\mm u}^\perp] \cdot {\mm u}^l)_e \\&= -\sum_{e\in {\mc E}} (c_e\, {\mm u}^l\cdot {\mm n}_e [{\mm u}],[{\mm u}])_e\,,
\end{aligned}
\end{equation}
where we used the fact that ${\mm u}^l$ is continuous and standard vector calculus identities. %,and the vector identity 
Similarly, taking ${\mm v}={\mm u}^s$ yields 
\begin{equation}
s({\mm u};{\mm u},{\mm u}^s) =  \sum_{e\in {\mc E}} (c_e\, {\mm u}^l\cdot {\mm n}_e [{\mm u}],[{\mm u}])_e\,.
\end{equation}
Due to the ${\mm L}^2$ orthogonality of ${\mm u}^l$ and ${\mm u}^s$, by taking ${\mm v}= {\mm u}^l$ and ${\mm v}={\mm u}^s$ in Equation~\eqref{eq:euler} we obtain the evolution equations for the large scale and small scale kinetic energy respectively, i.e.\ $E^l \coloneqq \|{\mm u}^l\|^2_\Omega/2$  and $E^s \coloneqq \|{\mm u}^s\|^2_\Omega/2$.  These are given by
\begin{equation} \label{eq:energies}
\left\{
\begin{array}{l}
\begin{aligned}
\displaystyle \der{E^l}{t} + a({\mm u};{\mm u},{\mm u}^l)&  -(p, \nabla\cdot {\mm u}^l)_\Omega\\& =  \sum_{e\in {\mc E}} (c_e\, {\mm u}^l\cdot {\mm n}_e [{\mm u}],[{\mm u}])_e\end{aligned}\\
\begin{aligned}
\displaystyle \der{E^s}{t}  + a({\mm u};{\mm u},{\mm u}^s)&  -(p, \nabla\cdot {\mm u}^s)_\Omega \\& = -\sum_{e\in {\mc E}} (c_e\, {\mm u}^l\cdot {\mm n}_e [{\mm u}],[{\mm u}])_e
\end{aligned}
\end{array}
\right.
\end{equation}
Then, if we set $c_e = \alpha\, {\mm u}^l \cdot n_e / (2|{\mm u}^l \cdot {\mm n}_e|)$ (by analogy to the usual choice in the flux form discretisation), we introduce an artificial energy transfer from small to large scales. This is because, for such a choice, $s({\mm u};{\mm u},{\mm u}^l)\leq 0$ whereas $s({\mm u};{\mm u},{\mm u}^s)\geq 0$. Moreover, the total energy is conserved since $s({\mm u};{\mm u},{\mm u}^l)+s({\mm u};{\mm u},{\mm u}^s)= s({\mm u};{\mm u},{\mm u})= 0$.

Finally, note that once the problem is reformulated as in Equation~\eqref{eq:energies}, the considerations of the previous section on the dissipation time scale apply without major changes. In particular, the dissipation time scale in the small scale equation is bounded from above by $\tau_h/\alpha$, with $\tau_h \coloneqq h/(\sup_\Omega\|{\mm u}^l\|_2)$.

\subsection{Resolved/unresolved interactions in the SUPG discretisation}\label{sec:supgscale}
Consider now the SUPG scheme. \col{We follow the variational multiscale (VMS) framework, and in particular \cite{codina11},} to give a simplified interpretation of the SUPG stabilisation as a subgrid scale  model. We assume that the exact vorticity $\omega \in {\mathcal V}_0$, where ${\mathcal V}_0$ is a sufficiently regular space, so that the following variational formulation,
\begin{equation}
\left\{
\begin{array}{l}
(\partial_t \omega, \varphi)_\Omega + ({\mm u}\cdot \nabla \omega, \varphi)_\Omega  = 0,  \\
{\mm u} = \nabla^\perp (\Delta^{-1}\omega),
\end{array}
\right.
\end{equation}
holds for all $\varphi \in {\mathcal V}_0 $.
Assume that ${\mathcal V}_0 = V_0 \oplus  V^u_0$, where $V^u_0$ contains the information from the unresolved scales and it is not necessarily orthogonal in $L^2$ to $V_0$. Then we can decompose $\omega = \bar{\omega} + \omega^u$, and by linearity ${\mm u}= \bar{{\mm u}}+ {\mm u}^u$. Using these decompositions, the vorticity equation becomes
\begin{equation}
\left\{
\begin{array}{l}
\begin{aligned}
(\partial_t \bar{\omega} +\partial_t \omega^u & + \bar{{\mm u}}\cdot \nabla \bar{\omega}+ \bar{{\mm u}}\cdot \nabla {\omega}^u, \bar{\varphi})_\Omega +\\ &( {\mm u}^u\cdot \nabla \bar{\omega} + {\mm u}^u\cdot \nabla {\omega}^u , \bar{\varphi})_\Omega  = 0, \end{aligned}\\
\begin{aligned}
(\partial_t \bar{\omega} + \partial_t {\omega}^u & + \bar{{\mm u}}\cdot \nabla \bar{\omega}+ \bar{{\mm u}}\cdot \nabla {\omega}^u, {\varphi}^u)_\Omega +\\&({{\mm u}}^u\cdot \nabla \bar{\omega} + {\mm u}^u\cdot \nabla {\omega}^u , {\varphi}^u)_\Omega  = 0,
\end{aligned}  \\
\end{array}
\right.
\end{equation}
for all $(\bar{\varphi},\varphi^u) \in V_0\times V_0^u$.
We now assume that the dynamics of the resolved and unresolved scales is driven only by the resolved velocity component. This assumption can be justified heuristically as in \cite{laval99}. If we also neglect the unsteady part of the unresolved component of the vorticity, we are left with the following system
\begin{equation}
\left\{
\begin{array}{l}
(\partial_t \bar{\omega}, \bar{\varphi})_\Omega + (\bar{{\mm u}}\cdot \nabla \bar{\omega}, \bar{\varphi})_\Omega - (\omega^u,\bar{{\mm u}}\cdot\nabla\bar{\varphi})_\Omega  = 0 \\
 (\partial_t \bar{\omega}, {\varphi}^u)_\Omega+(\bar{{\mm u}}\cdot \nabla \bar{\omega}+  \bar{{\mm u}}\cdot \nabla {\omega}^u, {\varphi}^u)_\Omega  = 0 \\
\end{array}
\right.
\end{equation}
for all $(\bar{\varphi},\varphi^u) \in V_0\times V_0^u$.
Finally, we approximate the unresolved scales equation by taking $(\bar{{\mm u}}\cdot \nabla \omega^u,\varphi^u)\approx \tau_s^{-1}(\omega^u,\varphi^u)$. Then, a possible solution of the unresolved scales equation is  $\omega^u = \tau_s R(\bar{\omega})$, with $R(\bar{\omega}) = \partial_t \bar{\omega} + \bar{{\mm u}} \cdot \nabla \bar{\omega}$. 
\col{The resulting scheme coincides with the SUPG scheme, with global intrinsic time scale $\tau_s$. In fact, this derivation shows that for our problem the VMS and SUPG methods are equivalent. This is due to the absence of the dissipation term in the equations studied here.}
\col{We remark that different versions of the above derivation have been used to devise similar algorithms.}
In  \cite{hughes95}, for example, a different analysis leads to the definition of a local time scale, which originates from the Green's function associated to the unresolved scales equation.  
\col{In \cite{codina02} the time dependent term in the unresolved scales equation is retained,  leading to a scheme in which the unresolved scales are dynamically tracked in time. These schemes still conserve energy as they share the same structure, however, in this paper we limit ourself to study the simpler SUPG scheme.}

It should be noted that the SUPG scheme models enstrophy transfer between different scales rather than energy transfer, as it was the case for the Lie derivative scheme. This is only possible due to the fact that we have an explicit equation for the vorticity evolution. Then, the diffusion introduced by the stabilisation term only affects enstrophy, and does not compromise energy conservation.  

\section{Numerical experiments}\label{sec:num}
%Energy/enstrophy tendency spectra

In this section we present some numerical results for the discretisations discussed in the previous sections. We start by reporting on their order of convergence using a manufactured solution test (Section~\ref{sec:order}). Then, we perform a forced turbulence simulation measuring energy and enstrophy tendencies in spectral space, in order to quantify the multiscale behaviour of the schemes (Section~\ref{sec:referencet} and \ref{sec:finitet}). We use as a benchmark test case the forced turbulence test performed in \cite{thuburn14}, and we can compare our results to their reference solution and to the other schemes tested therein. \col{Note that as in \cite{thuburn14} for all the tests of this section we regard the equations of motion as non-dimensional, since they can always be brought in this form via an appropriate time rescaling.}
All computations were performed using the Firedrake software suite \citep{Rathgeber2015}, which allows for symbolic implementation of finite element problems of mixed
type. We also used additional resources for the implementation from the following references \citep{petsc-user-ref,petsc-efficient,Dalcin2011,Chaco95,MUMPS01,MUMPS02}.

\subsection{Order of convergence}\label{sec:order}
The method of manufactured solutions provides a useful way to investigate
convergence rates when analytic solutions are not available. 
\col{This consists in choosing an arbitrary function and adding an appropriate forcing to the system of equations we want to solve, so that this function is an exact solution of the modified system.}
\col{We let $\Omega = [0,1]\times [0,1]$, with Cartesian coordinates $(x,y)$}, and pick as manufactured solution the one generated by the following stream function,
\begin{equation}\label{eq:manufactured}
{\psi}(t,x,y) = \sin(\pi x)\sin(\pi y)\sin(\pi y-t)e^{-2t/\sigma} \, ,
\end{equation}
with $\sigma=100$, $t\in [0,1]$ and boundary conditions $\psi = 0$ on $\partial \Omega$. The $L^2$ errors in velocity and vorticity are reported in Table~\ref{tab:centup} and \ref{tab:centupo} respectively. For the Lie derivative scheme we set $\alpha=1$,  whereas for the SUPG scheme the stabilisation coefficient is given by Equation~\eqref{eq:tau} with $\beta=1$. Both schemes are integrated in time using the implicit midpoint integration rule with $\Delta t = 10^{-3}$. Note that for the Lie derivative scheme the vorticity needs to be computed weakly at each time step since it is not used directly by the algorithm. The schemes behave similarly in terms of velocity, where $r+1$ order of convergence is observed. However, the SUPG scheme yields lower error and convergence rate greater than $r+1$ for the vorticity, whereas the Lie derivative scheme only gives approximately order $r$ convergence.

\begin{table}[t!]
\centering
\begin{tabular}{c c c c c c}
\toprule
 \multirow{2}{*}{$r$}& \multirow{2}{*}{$h$} & \multicolumn{2}{c}{SUPG}& \multicolumn{2}{c}{Lie derivative}\\
& & error
& order& error
& order\\
\midrule
\multirow{4}{*}{1}& 7.40e-1 & 2.68e-2 &-- & 4.12e-2& -- \\
& 3.70e-1 & 
 6.69e-3& 2.00 &
  5.38e-3
   & 2.19\\
& 2.47-1 & 
2.97e-3
& 2.00 & 9.00e-3
& 2.10\\
& 1.85e-1 & 1.67e-3& 2.00&   3.84e-3& 2.09\\
\midrule
\multirow{4}{*}{2}& 7.40e-1 & 1.37e-3 &-- & 3.03e-3& -- \\
& 3.70e-1 & 
 1.71e-4& 3.01 &
   3.20e-4 & 3.25 \\
& 2.47-1 & 
  5.05e-5
& 3.00 &
   8.94e-5
& 3.14\\
& 1.85e-1 &
  2.13e-5& 3.00& 3.63e-5& 3.13\\
\bottomrule
\end{tabular}
\caption{Comparison between the SUPG and Lie derivative scheme in terms of the error $\|{{\mm u}}-{\mm u}_h\|_{\Omega}$ and the order of convergence, for the solution in Equation~\ref{eq:manufactured}, and for $\alpha=\beta=1$.}
\label{tab:centup}
\end{table}

\begin{table}[t!]
\centering
\begin{tabular}{c c c c c c}
\toprule
 \multirow{2}{*}{$r$}& \multirow{2}{*}{$h$} & \multicolumn{2}{c}{SUPG}& \multicolumn{2}{c}{Lie derivative}\\
& & error
& order& error
& order\\
\midrule
 \multirow{4}{*}{1}& 7.40e-1 & 1.04e-1 &-- & 1.58& -- \\
& 3.70e-1 & 
   2.16e-2
   & 2.27 &
   8.05e-1
   & 0.97\\
& 2.47-1 & 
8.62e-3
& 2.26 & 5.51e-1
& 0.93\\
& 1.85e-1 & 
    4.48e-3&  2.28&   4.21e-1& 0.94\\
    \midrule
     \multirow{4}{*}{2}& 7.40e-1 & 3.50e-3&-- &  1.34e-1& -- \\
    & 3.70e-1 &   2.79e-4
       & 3.65 &
        2.69e-2
             & 2.32\\
    & 2.47-1 & 
    6.74e-5
    & 3.50 &  1.10e-2    
    & 2.21\\
    & 1.85e-1 &   2.51e-5 & 3.43&   5.79e-3& 2.22\\
\bottomrule
\end{tabular}
\caption{Comparison between the SUPG and Lie derivative scheme in terms of the error $\|{{\omega}}-{\omega}_h\|_{\Omega}$ and the order of convergence,  for the solution in Equation~\ref{eq:manufactured}, and for $\alpha=\beta=1$.}
\label{tab:centupo}
\end{table}

\subsection{Reference turbulence test}\label{sec:referencet}
We now describe the reference turbulent solution considered in \cite{thuburn14}. We start by modifying the governing equation by adding a forcing and a dissipation term as follows:
\begin{equation}\label{eq:eulerf}
\displaystyle \partial_t{\omega} + u \cdot \nabla \omega   =  f -\omega/\tau\, ,
\end{equation}
on a periodic unit square domain $\Omega$, with Cartesian coordinates $(x,y)$. The fixed forcing is $f = 0.1\sin(32\pi x)$ and acts on the wave number $k =16$. The dissipation term is the standard Ekman drag with time scale $\tau = 100$, which does not introduce selective decay \citep{majda06}. The initial conditions are
\begin{equation}\label{eq:ic}
\begin{array}{ll}
\omega|_{t=0} =&\sin(8\pi x)\sin(8\pi y)\\& + 0.4 \cos(6\pi x)\cos(6\pi y) \\&+ 0.3 \cos(10\pi x)\cos(4\pi y) \\&
+0.01\sin(2\pi y) +0.02 \sin(2\pi x)\,,
\end{array}
\end{equation}
after which the flow rapidly becomes chaotic.

\begin{figure*}
\centering
\includegraphics[scale=0.9]{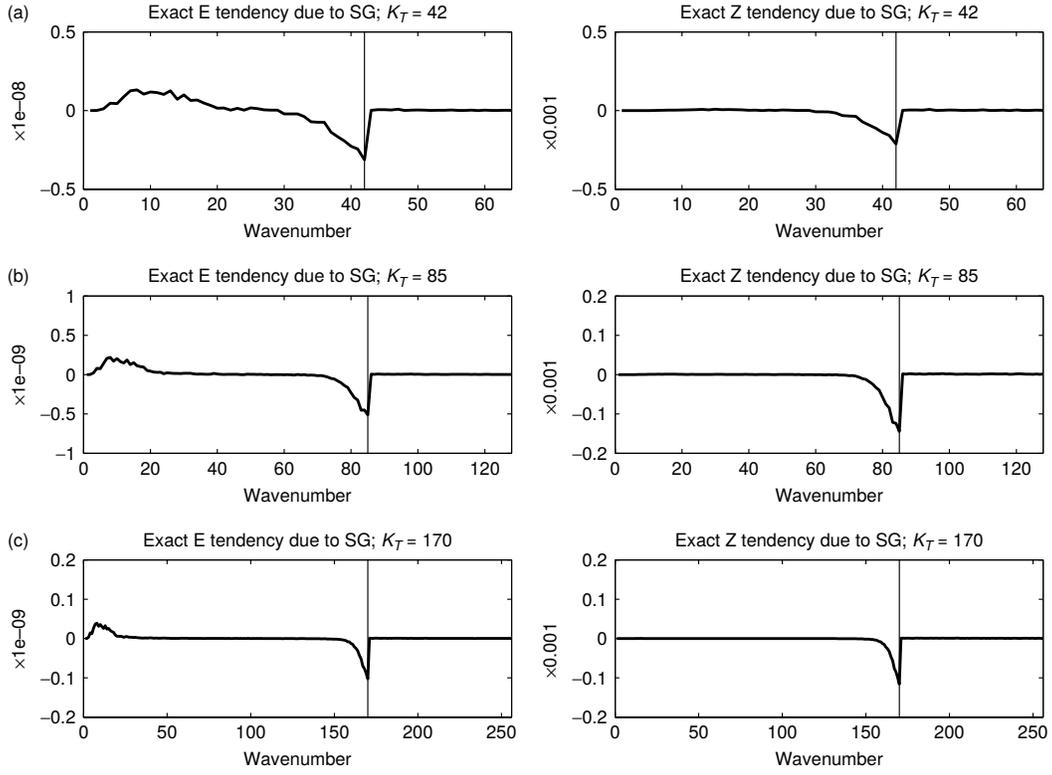}
\caption{Energy and enstrophy tendencies for the reference solution, with $k_T= 42$ (top), $k_T =85$ (middle), $k_T =170$ (bottom) (Reproduced from \citet{thuburn14}). }\label{fig:thuburn_Fig1}
\end{figure*}

\cite{thuburn14} analysed the stationary state in terms of energy and enstrophy tendencies in spectral space, with the aim of quantifying the effect of small scales (above a certain threshold) on large scales. In order to compute the reference values, they advanced the solution of Equation~\eqref{eq:eulerf} until $t=200$ using a spectral method with high resolution, computing the stream function $\psi$ and the vorticity $\omega$.  Then, for $t\in[200,210]$ they computed  the Jacobian $J = \nabla^\perp \psi \cdot \nabla\omega$ and then the energy and enstrophy time derivative at each time step, using the following formulas
\begin{equation}
\dot{E}(k_x,k_y) = \mr{Re} \left\{ \frac{\hat{\psi}^* \hat{J}}{ N^4}  \right\}, \quad
\dot{Z}(k_x,k_y) = \mr{Re} \left\{ \frac{\hat{\omega}^* \hat{J}}{ N^4}  \right\},
\end{equation}
\col{where $N$ is the number of grid points}, $\hat{\psi}$ and $\hat{J}$ are the Fourier transforms of $\psi$ and $J$ (truncated at the maximum retained wavenumber), and the superscript $*$ denotes complex conjugation. The results were integrated over angle in spectral space to give $\dot{E}(k)$ and $\dot{Z}(k)$, where $k\coloneqq \sqrt{k_x^2+k_y^2}$. The same procedure was repeated truncating all the Fourier transforms with $k\leq k_T$, for several values of $k_T$, yielding $\dot{E}_T(k)$ and $\dot{Z}_T(k)$. Finally, the energy and enstrophy tendencies due to scales $k>k_T$ were computed as follows:
\begin{equation}\label{eq:SG}
\dot{E}_{SG}(k) = \dot{E}(k)-\dot{E}_T(k),~~~ \dot{Z}_{SG}(k) = \dot{Z}(k)-\dot{Z}_T({k})  . \end{equation}
The results of such calculations are illustrated in Figure~\ref{fig:thuburn_Fig1}. The main features are a net non-local energy transfer from small to large scales and enstrophy dissipation concentrated close to the truncation wavenumber. Reproducing this behaviour on the resolved scale of a numerical simulation is very challenging. %(Performances of the schemes in \cite{thuburn14})[...] 
In the following, we review the behaviour of the schemes presented in the previous section.

\subsection{Finite element turbulence tests}\label{sec:finitet}
We now repeat the experiment described in the last section, but using the methods defined in Section~\ref{sec:methods}. In practice, we substitute $\dot{E}(k)$ and $\dot{Z}(k)$ in Equation~\eqref{eq:SG} with those computed using the numerical Jacobian implied by the discretisations.  Instead of restarting each simulation for the initial condition specified by a DNS solution at $t=200$ as in \cite{thuburn14}, we start monitoring the several schemes at a later time $t=300$, so that each simulation can be considered to have reached its own statistical equilibrium. \col{In particular, we verified this by checking the time evolution of the energy moving average.} For all simulations we used a structured right triangular grid and we denote by $h$ the horizontal and vertical spacing of the nodes. Moreover, they are all integrated in time using the implicit midpoint integration rule with $\Delta t = 2\cdot 10^{-2}$ until $t=350$.

\col{We compare the results for the flux form scheme, the Lie derivative scheme and the SUPG scheme when using the same function spaces. The results for the flux form scheme are mainly reported for comparison with the Lie derivative scheme, since these schemes share a similar structure and the same convergence rate. In all tests, the stabilisation factor for the SUPG scheme is given by Equation~\eqref{eq:tau} with $h_T=h/\sqrt{2}$, as suggested by \cite{codina92} for a symmetrical grid of right triangles.}

\begin{figure*}[!h]
\centering
\includegraphics[scale=1.0]{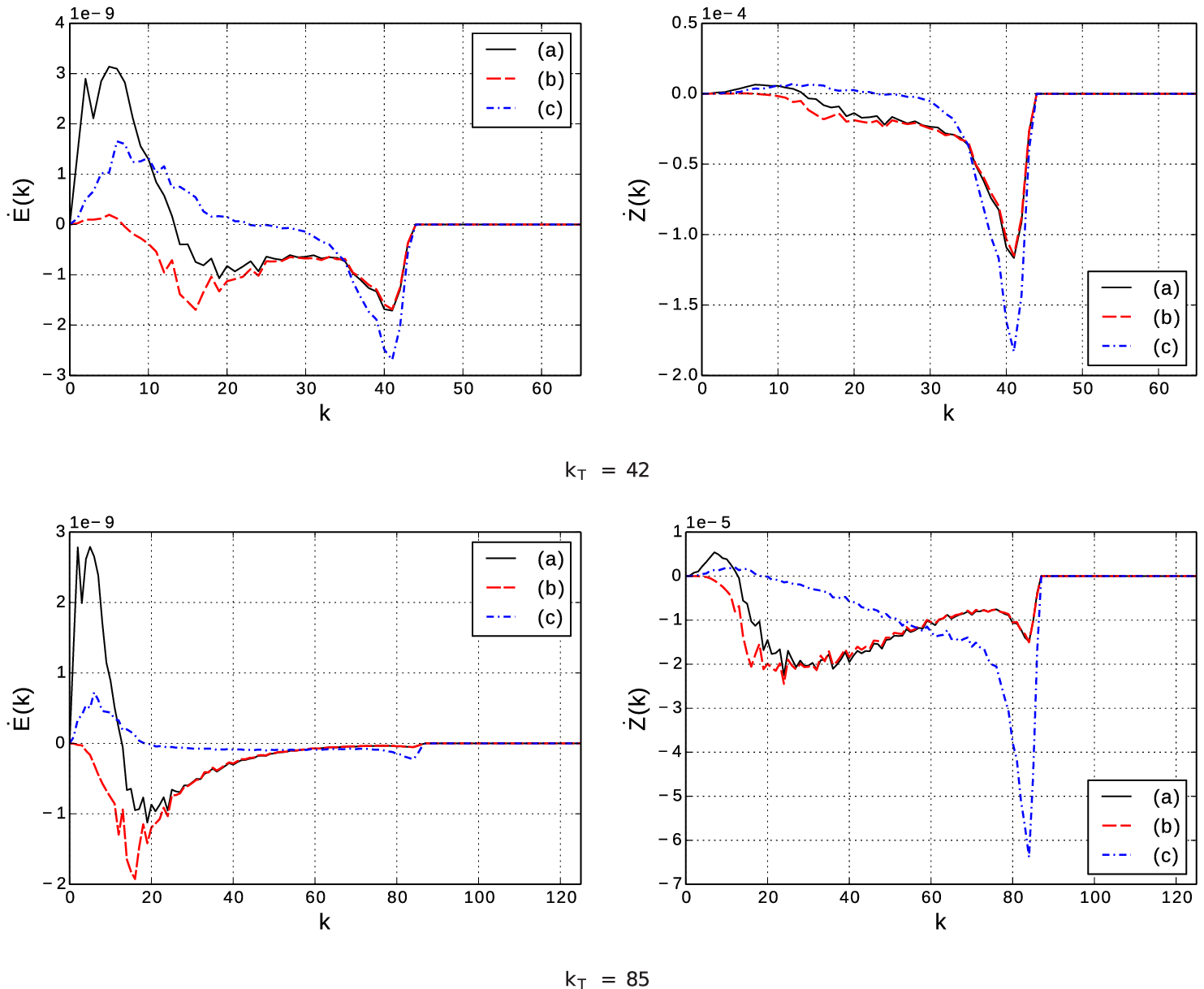}
\caption{Energy and enstrophy tendencies for the Lie derivative (a), the flux form (b) and the SUPG (c)  scheme (with $h=1/128$, $r=1$, $\alpha=1$,$\beta=1$) on selected spectral intervals $k\leq k_T$.} \label{fig:EZ0}
\end{figure*}

In Figure~\ref{fig:EZ0} results are shown for $r=1$, $h= 1/128$, $\alpha=1$ and $\beta=1$. The mesh resolution is such that the maximum resolved wavenumber (computed using the grid associated to the stream function degrees of freedom) is approximately $k=85$. For this setting we show the effect of the numerical scheme on the integration for the cut-off wavenumbers tested in \cite{thuburn14}. 

We start by analysing the energy behaviour (left column in Figure~\ref{fig:EZ0}). At the lowest cut-off wavenumber $k_T= 42$, the positive energy peak is overestimated for the Lie derivative scheme but it is completely absent when using the flux form discretisation. Both schemes are too dissipative at intermediate scales, with the flux form discretisation being more dissipative at the forcing scale. For $k_T= 85$ the behaviour of the two schemes is substantially different. The flux form scheme appears to subtract energy over the whole spectrum, and especially at the most energetic wavenumbers (since energy is inserted at $k=16$). The Lie derivative scheme also subtracts energy over a larger wavenumber interval compared to the reference test in Figure~\ref{fig:thuburn_Fig1}, but energy is reinserted at low frequencies, as expected from the discussion in the previous section. However, note that the positive peak is considerably larger than the one in the reference solution. 

The SUPG scheme shows a much better agreement to the reference test, at equal resolution. Note, in particular, that for $k_T = 85$ we see an energy dissipation peak close to the cut-off wavenumber, which was absent for the other schemes tested. This is probably due to the fact that the vorticity is explicitly used in the computation of the Jacobian, leading to a lower error in the computation of the energy and enstrophy tendencies. Such a conjecture is corroborated by the order of convergence results given in Section~\ref{sec:order}. Nonetheless, most of the dissipation is still acting on lower wavenumbers, ranging almost down to the forcing scale.

As for the enstrophy behaviour (right column in Figure~\ref{fig:EZ0}), the Lie derivative and the flux form schemes show similar behaviour for $k_T=42$, where enstrophy is subtracted mainly at wavenumbers close to the cut-off wavenumber as in Figure~\ref{fig:thuburn_Fig1}. For $k_T=85$ however, even if the negative peak at the highest wavenumber is still represented, this is much smaller than in the reference solution, and enstrophy appears to be dissipated comparatively more over the rest of the spectrum for both schemes. \col{Again, the SUPG scheme produces a much better agreement to the reference solution, and it is able to capture the sharp dissipation enstrophy peak also for $k_T =85$.}

\begin{figure*}[!h]
\centering
\includegraphics[scale=1.0]{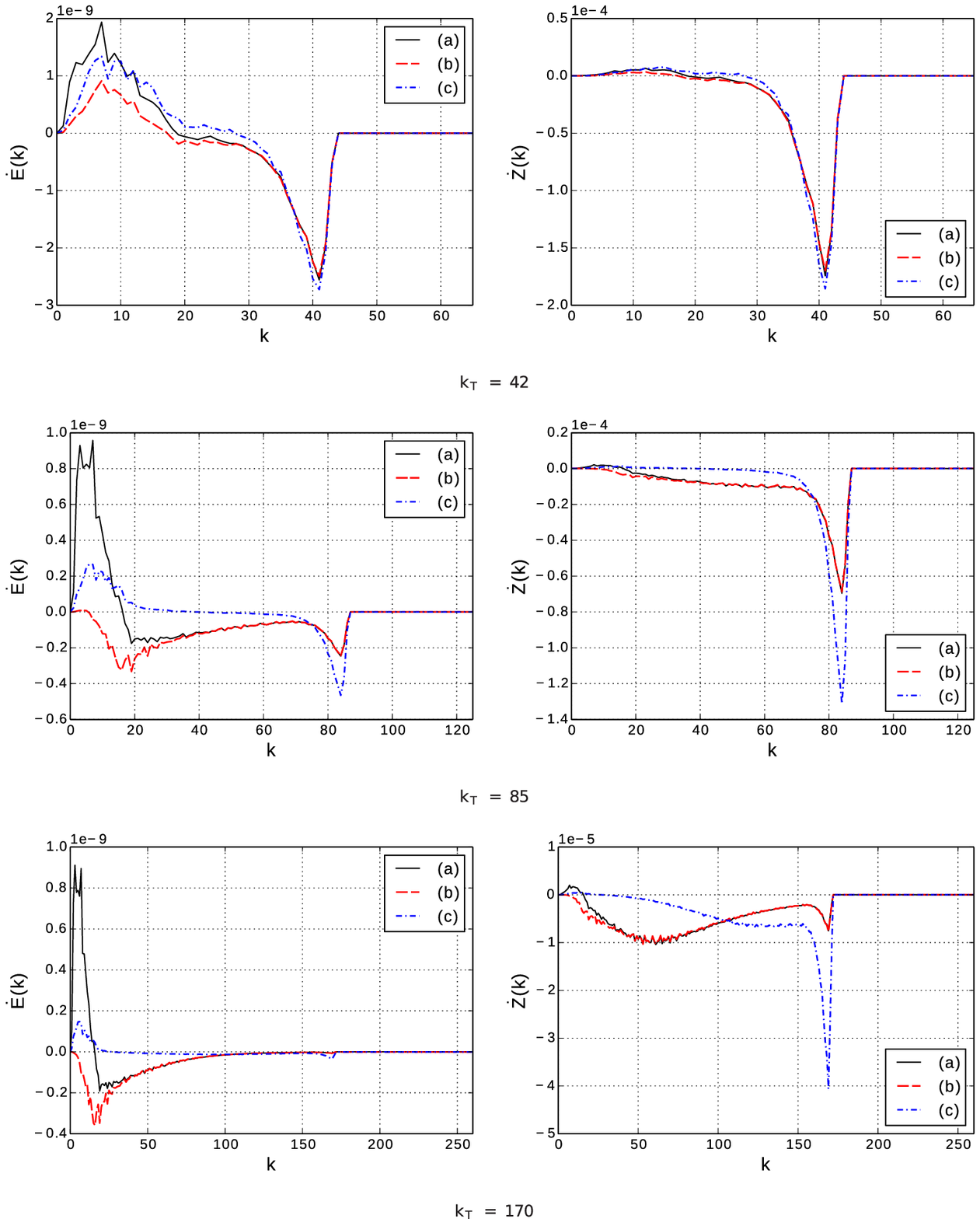}
\caption{Energy and enstrophy tendencies for the Lie derivative (a), the flux form (b) and the SUPG (c)  scheme (with $h=1/256$, $r=1$, $\alpha=1$,$\beta=1$) on selected spectral intervals $k\leq k_T$.} \label{fig:EZ1}
\end{figure*}

In Figure~\ref{fig:EZ1} the same test is repeated at higher resolution $h=1/256$ and same polynomial degree $r=1$, with maximum retained wavenumber approximately equal to $k=170$. For $k_T=42$,
the resolved scales dominate the energy and enstrophy tendencies for all the schemes, so that we can clearly see the non-local energy backscatter correctly represented. However, the energy conserving schemes show again a larger positive peak at low frequencies. For $k_T=85$, we can see that the Lie derivative scheme produces a positive energy peak, which is closer to its reference value, if compared to the coarser resolution test in Figure~\ref{fig:EZ0}. The peak is absent in the results for flux form discretisation. Nonetheless, the plot for $k_T=170$ confirms that the energy transfer is still strongly local as dissipation is mostly acting close to the forcing scale. \col{The SUPG scheme, on the other hand, does not overestimate the energy peak for $k_T = 85$, and produces a much more nonlocal energy transfer, as it can be verified from the plots for $k_T =170$.}

\begin{figure*}[!h]
\centering
\includegraphics[scale=1.0]{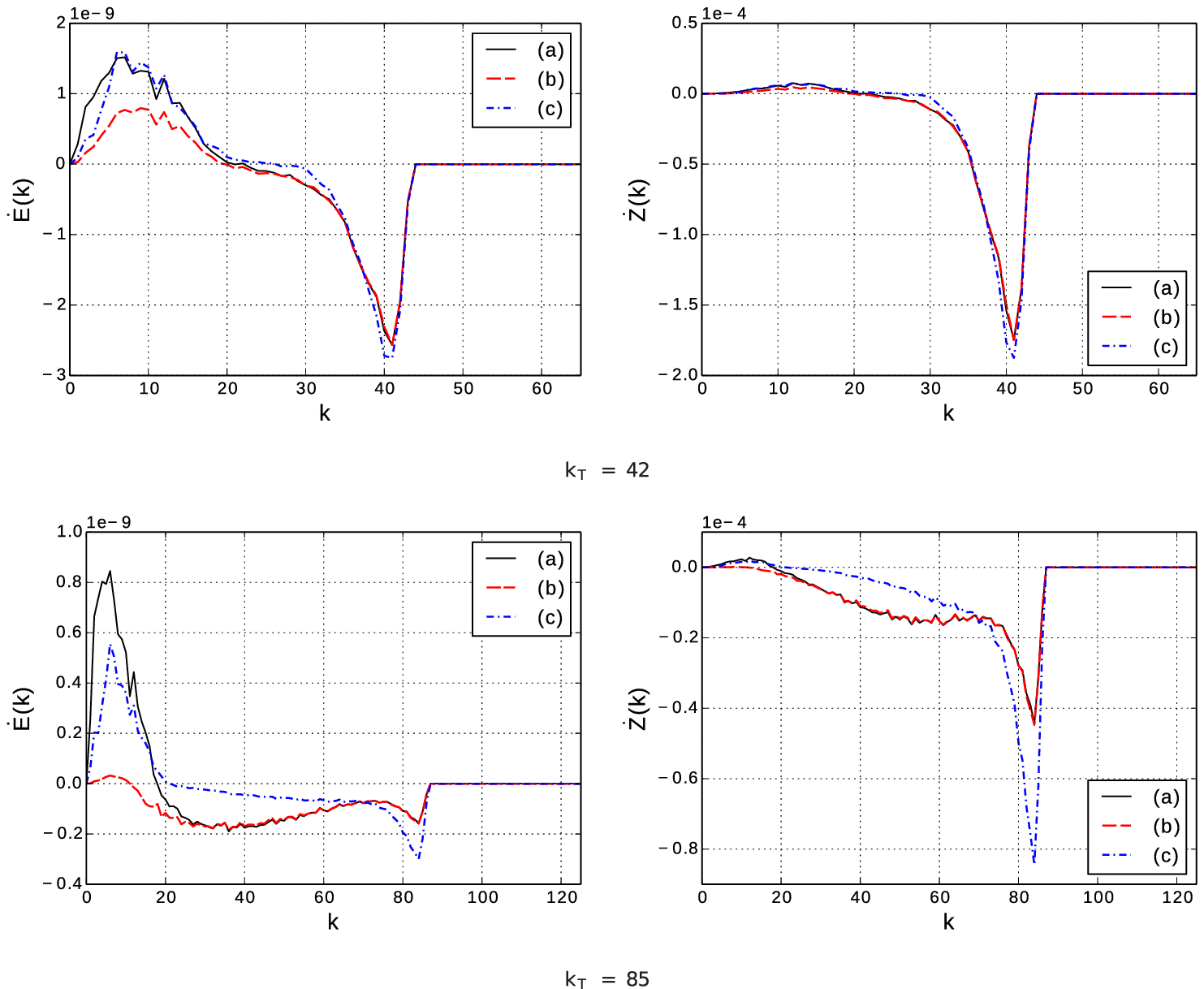}
\caption{Energy and enstrophy tendencies for the Lie derivative (a), the flux form (b) and the SUPG (c)  scheme (with $h=1/86$, $r=2$, $\alpha=1$,$\beta=1$) on selected spectral intervals $k\leq k_T$.} \label{fig:EZ2}
\end{figure*}

Next, in Figure~\ref{fig:EZ2}, we performed another test with large mesh size, i.e.\ we se set $h=1/85$, but with higher polynomial degree $r=2$. For this case the number of degrees of freedom is the same as the test case in Figure~\ref{fig:EZ0}. The behaviour, however, is substantially better for the Lie derivative scheme, with energy being dissipated over larger wavenumbers and the positive peak being comparable to the higher resolution tests in Figure~\ref{fig:EZ1}.
The improvement is more modest for the SUPG scheme and does not compare to the increased accuracy achieved with mesh refinement.
% % % %

%In Figure~\ref{fig:EZS2} we show the effect of mesh refinement by setting $h=256$ and $r=1$. The results are again significantly better than the ones obtained with the Lie derivative scheme. In Figure~\ref{fig:EZS3} we use higher order polynomials, i.e.\ we set $h=85$ and $r=3$, so that the highest resolved wavenumber is the same as for the results in Figure~\ref{fig:EZS}. 

\begin{figure*}[!t]
\centering
\subfigure[Flux form]{\includegraphics[clip =true, trim = {308 60 308 60}, scale=0.35]{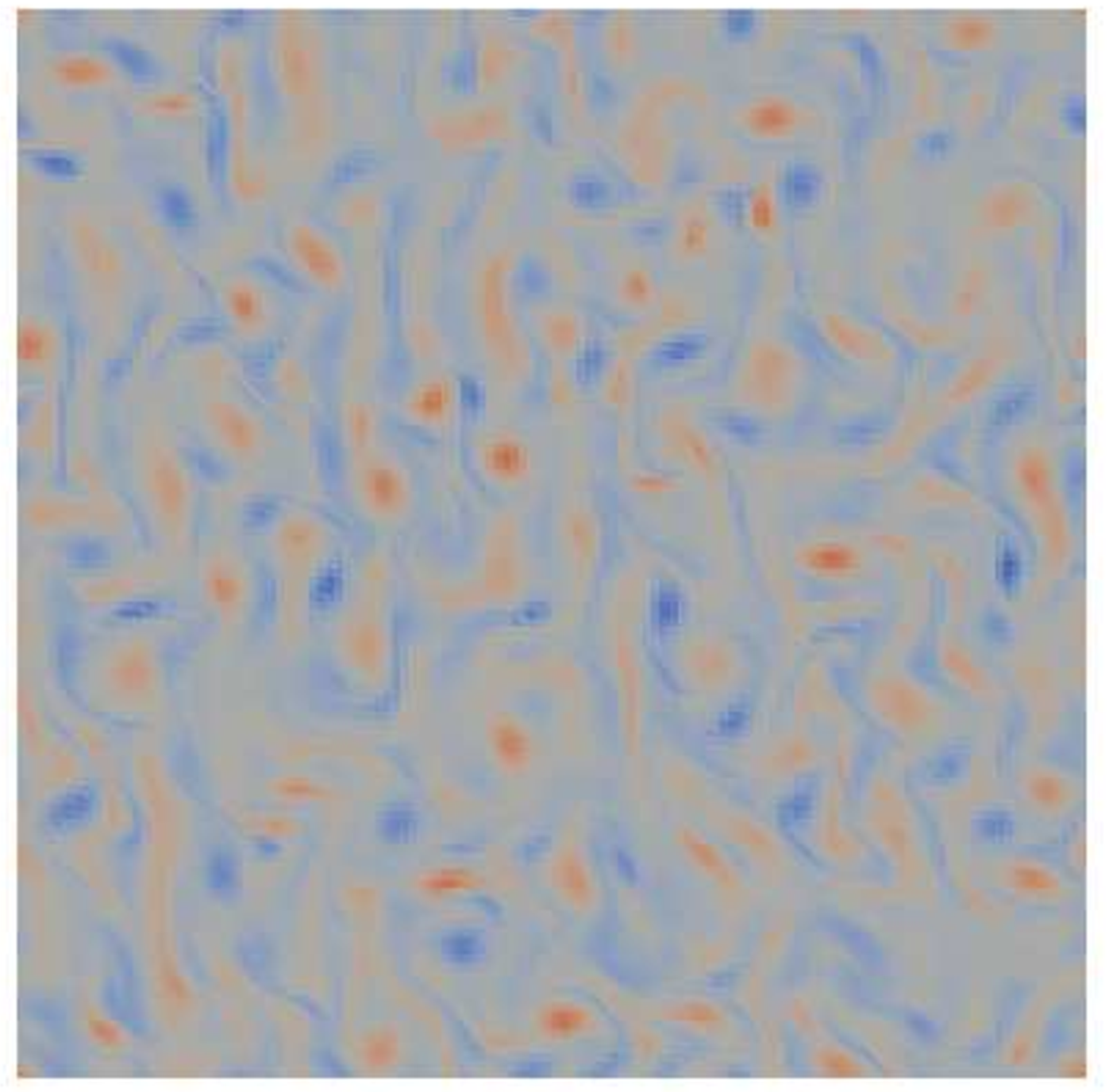}}
\subfigure[Lie derivative]{\includegraphics[clip =true, trim = {308 60 308 60}, scale=0.35]{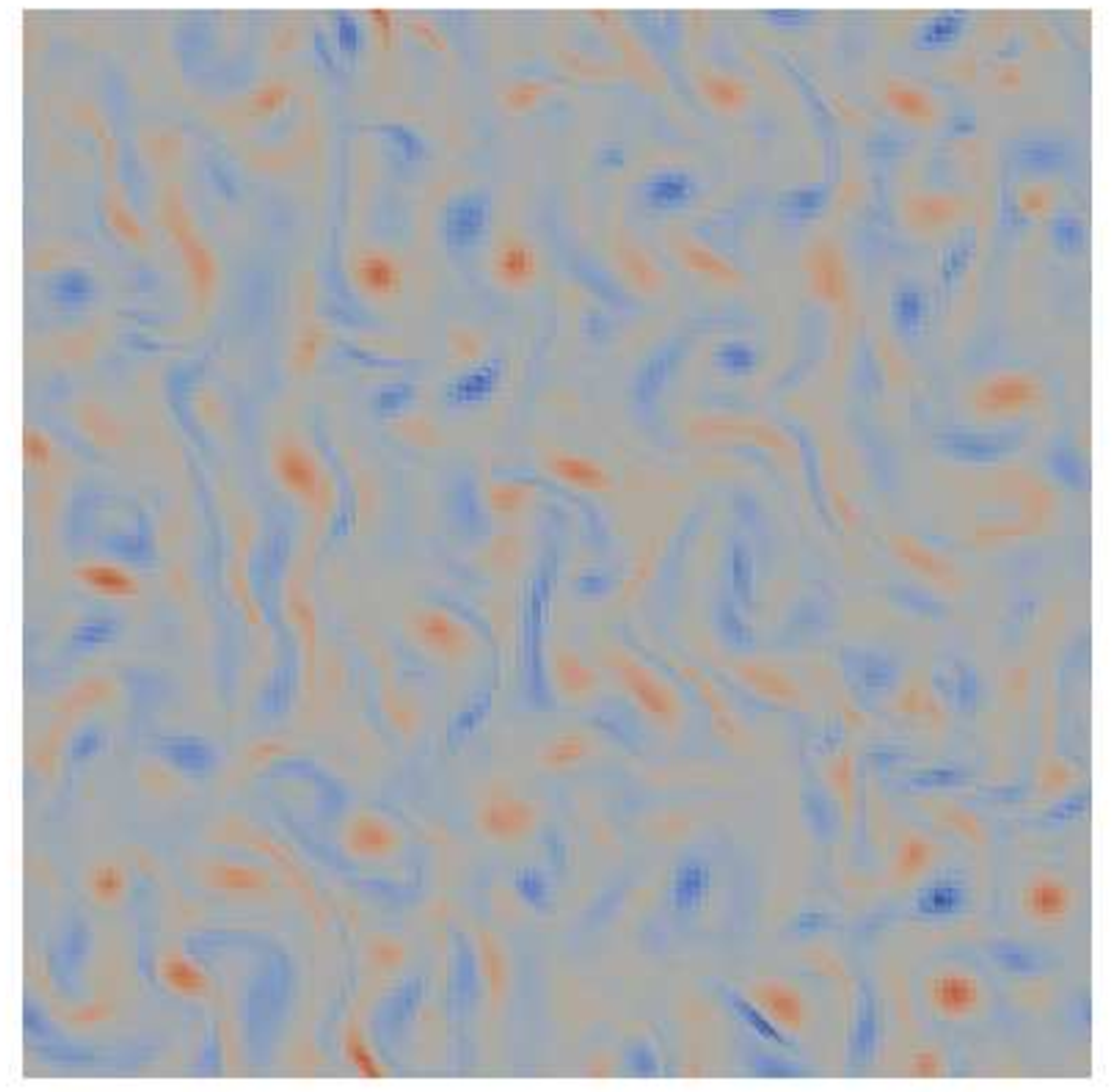}}
\subfigure[SUPG]{\includegraphics[clip =true, trim = {308 60 308 60}, scale=0.35]{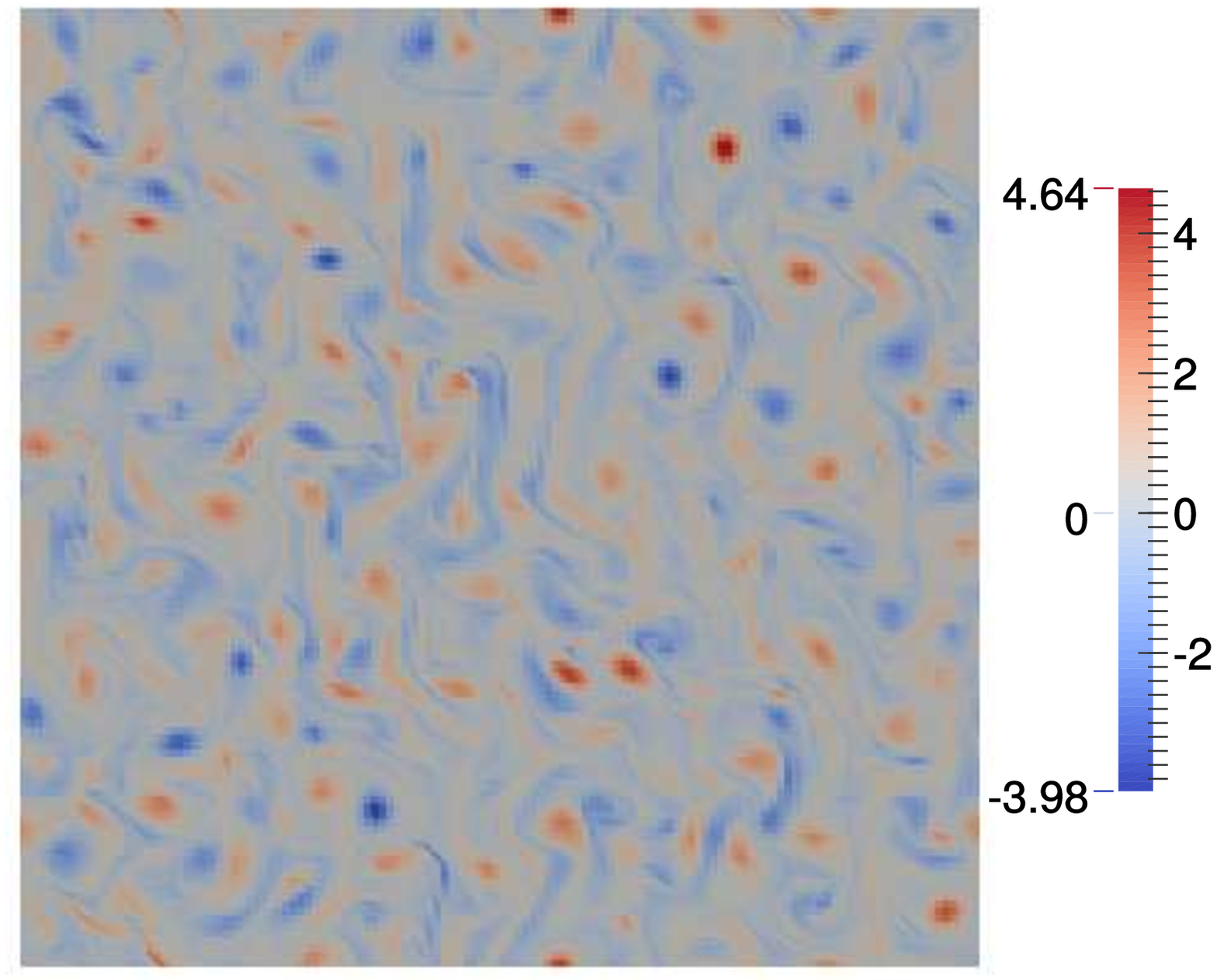}}
\subfigure{\includegraphics[clip =true, trim = {702 60 160 60}, scale=0.33]{Vort_128S.eps}}
\caption{Vorticity field at $t=300$, for the tested schemes (with $h=1/128$, $r=1$, $\alpha=1$,$\beta=1$).} \label{fig:Vort1}
\end{figure*}

\begin{figure*}[!t]
\centering
\subfigure[Flux form]{\includegraphics[clip =true, trim = {308 60 308 60}, scale=0.35]{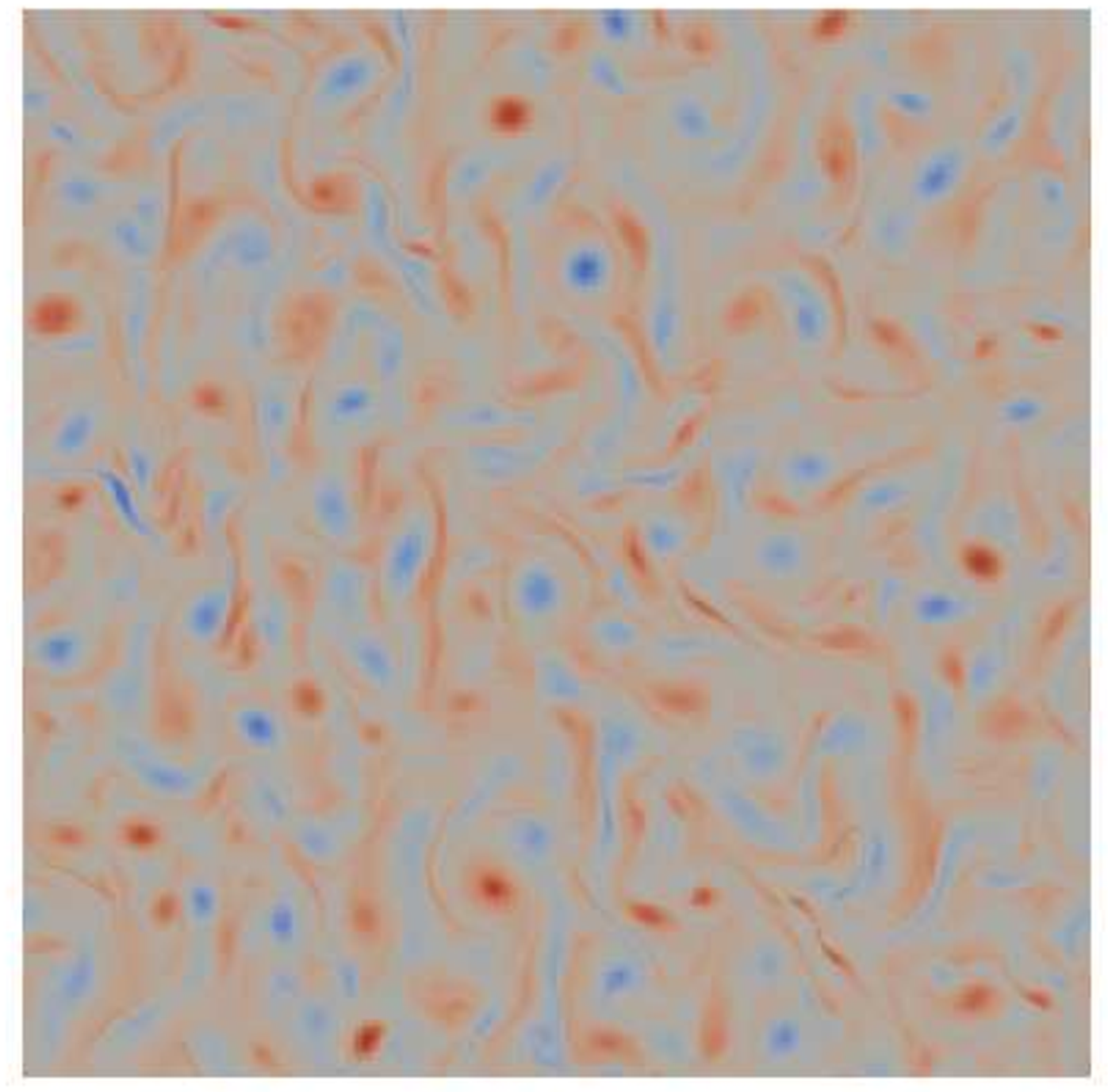}}
\subfigure[Lie derivative]{\includegraphics[clip =true, trim = {308 60 308 60}, scale=0.35]{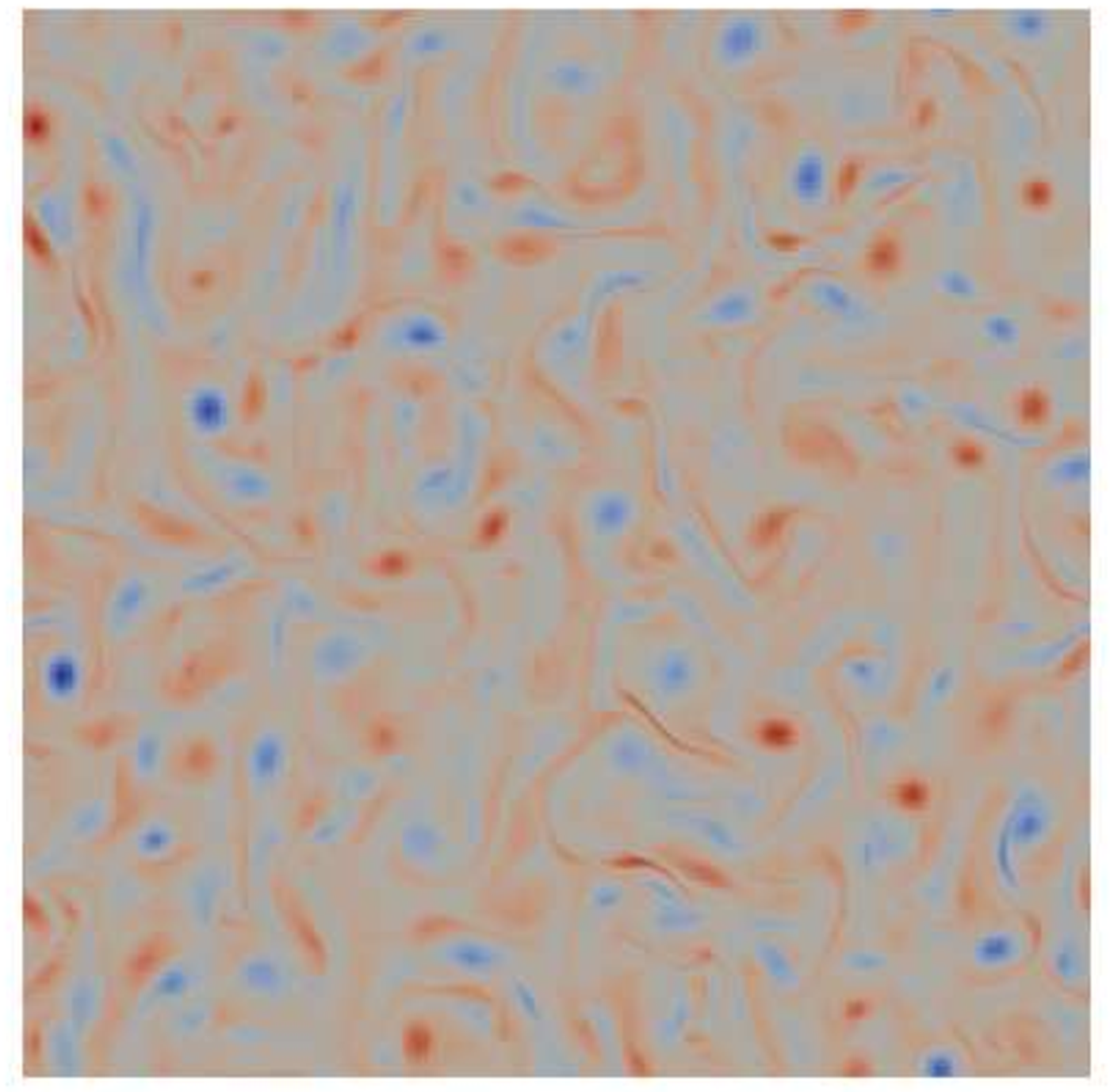}}
\subfigure[SUPG]{\includegraphics[clip =true, trim = {308 60 308 60}, scale=0.35]{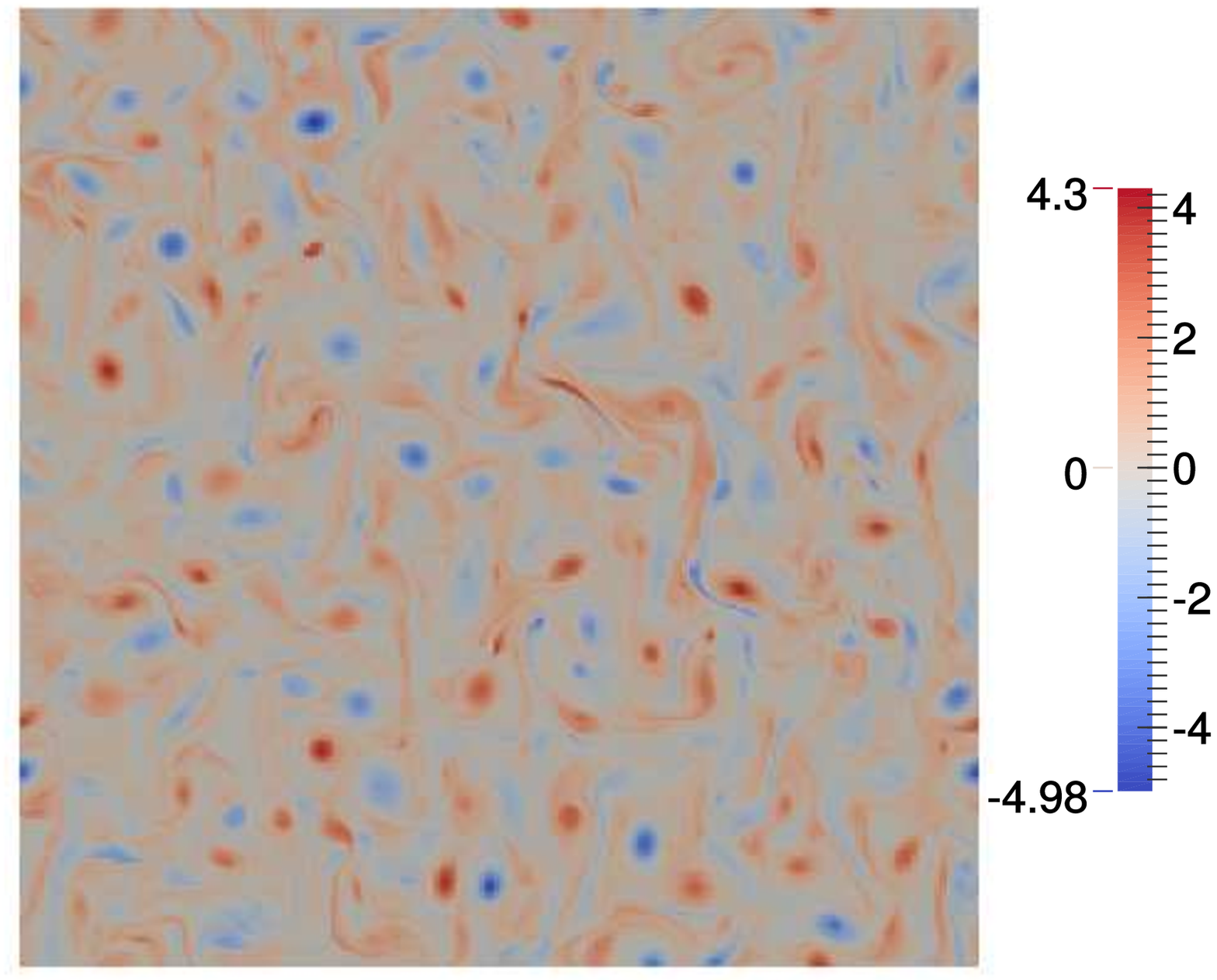}}
\subfigure{\includegraphics[clip =true, trim = {702 60 160 60}, scale=0.33]{Vort_256S.eps}}
\caption{Vorticity field at $t=300$, for the tested schemes (with $h=1/256$, $r=1$, $\alpha=1$,$\beta=1$).} \label{fig:Vort2}
\end{figure*}

\begin{figure*}[!t]
\centering
\subfigure[Flux form]{\includegraphics[clip =true, trim = {308 60 308 60}, scale=0.35]{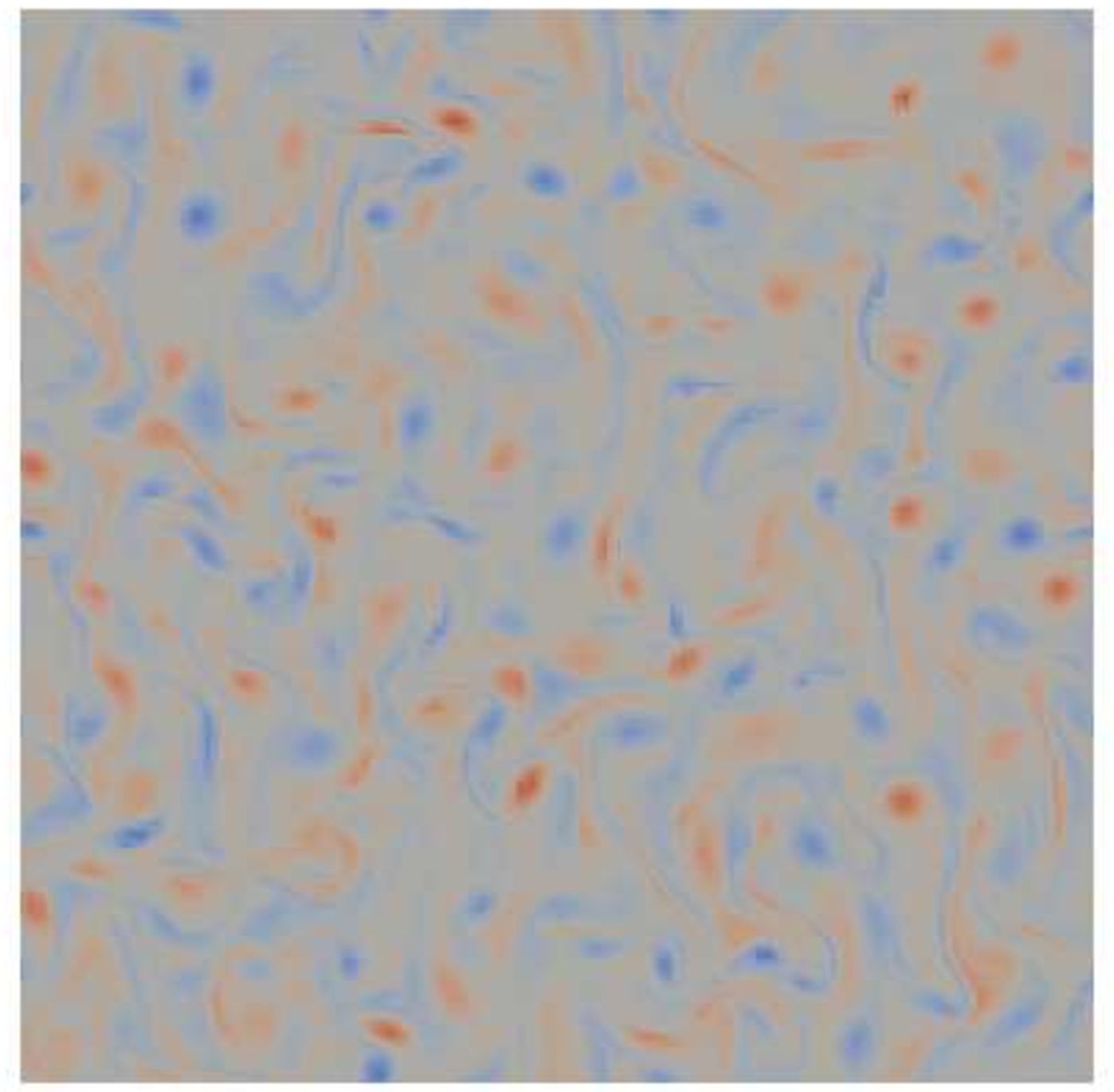}}
\subfigure[Lie derivative]{\includegraphics[clip =true, trim = {308 60 308 60}, scale=0.35]{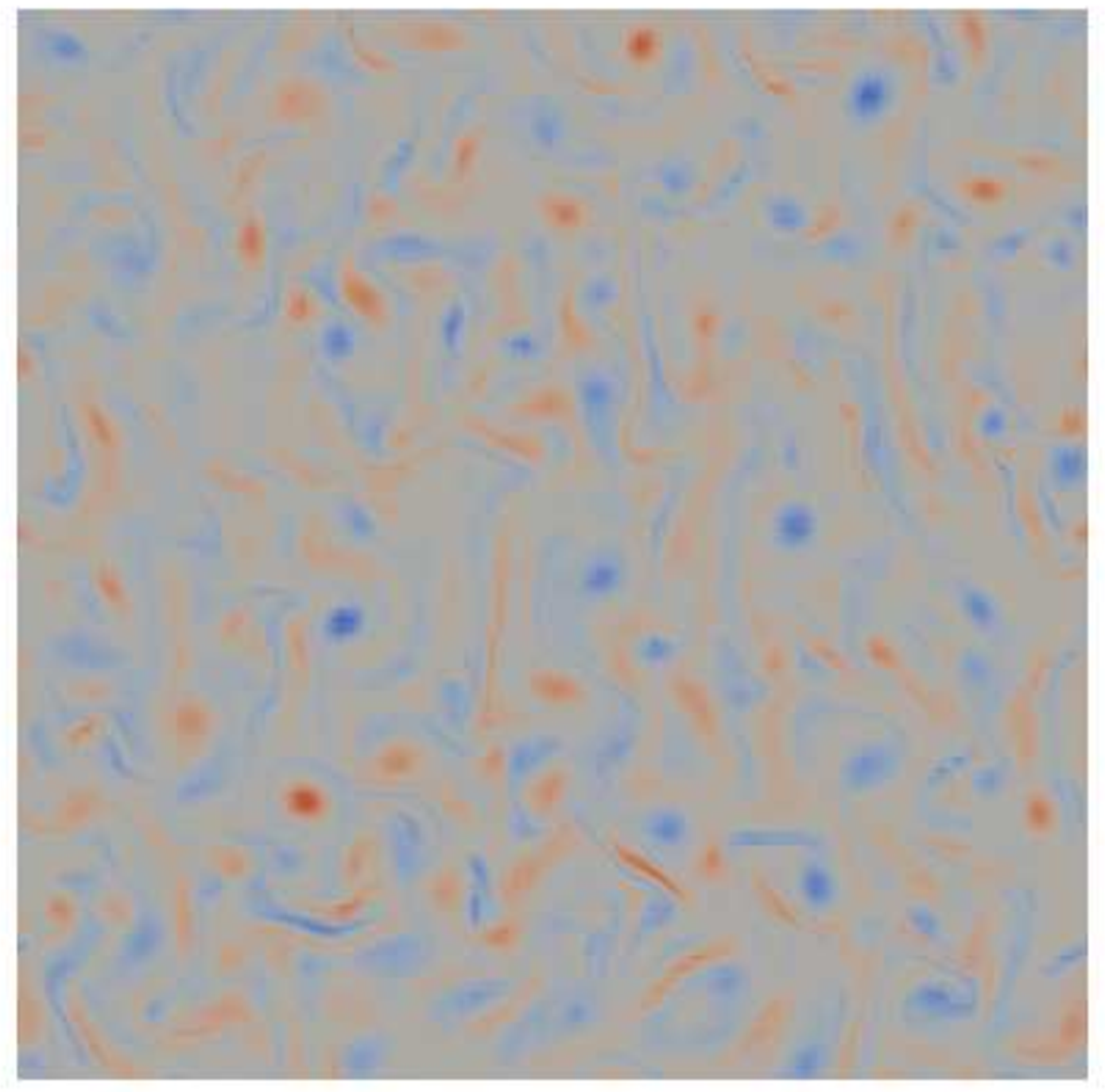}}
\subfigure[SUPG]{\includegraphics[clip =true, trim = {308 60 308 60}, scale=0.35]{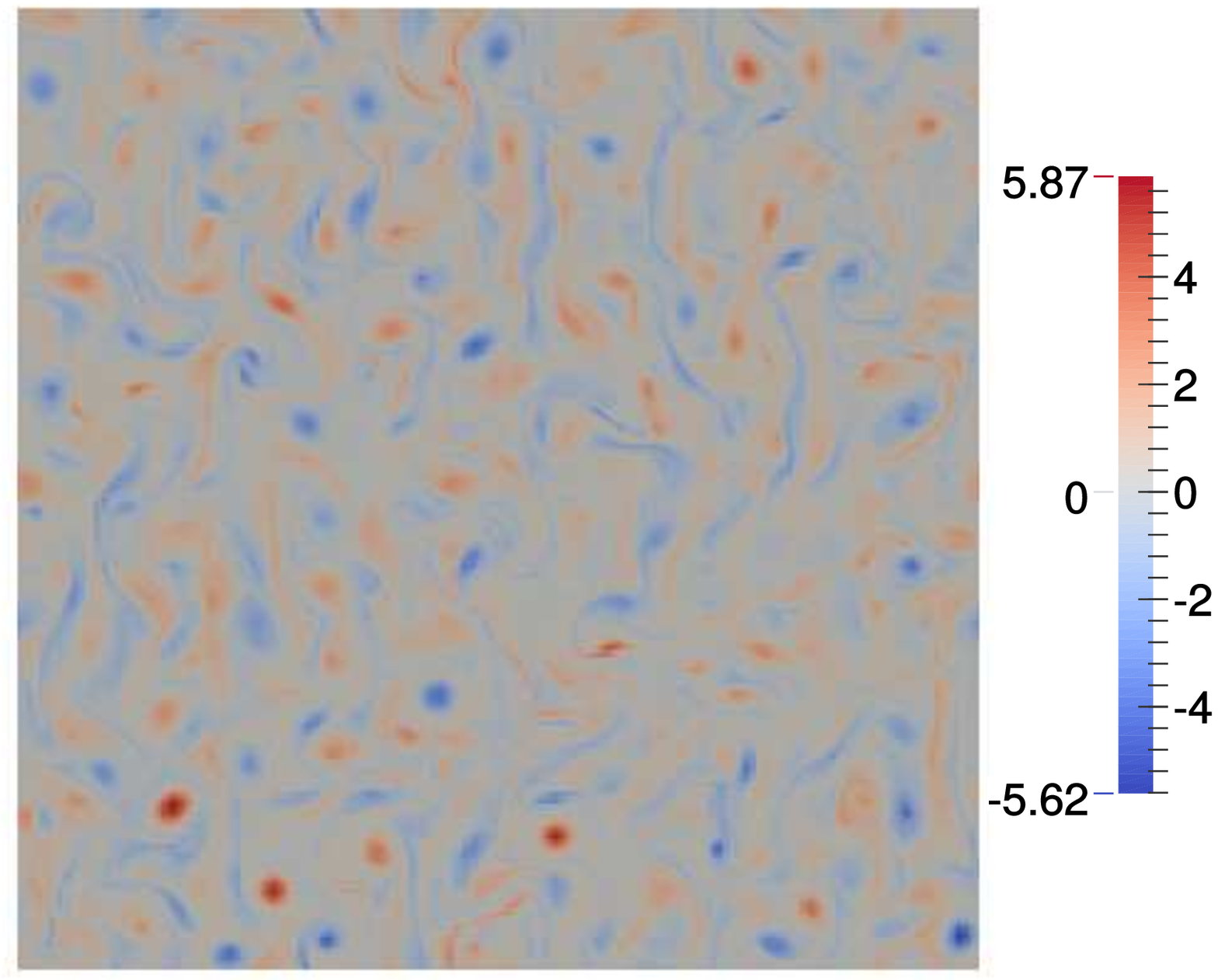}}
\subfigure{\includegraphics[clip =true, trim = {702 60 160 60}, scale=0.33]{Vort_86S.eps}}
\caption{Vorticity field at $t=300$, for the tested schemes (with $h=1/86$, $r=2$, $\alpha=1$, $\beta=1$).} \label{fig:Vort3}
\end{figure*}

\col{It should be noted that for both the Lie derivative and the SUPG schemes the upscale energy transfer is often accompanied by a spurious enstrophy injection at large scales. This effect is particularly pronounced in the low resolution test in Figure~\ref{fig:EZ0}.
In order to obtain an indication of the consequences of such a behaviour we compare the vorticity fields for the different schemes in physical space at fixed time $t=300$, see Figures~\ref{fig:Vort1}, \ref{fig:Vort2} and \ref{fig:Vort3}. We observe that, especially at low resolution ($h=1/128$, $r=1$), the vorticity field for the Lie derivative scheme appears smoother than the one obtained using the flux form scheme. Their behaviour however is comparable, whereas the SUPG scheme provides a more accurate representation of the vorticity field, with smaller and more dense vortices.
This reflects the higher accuracy of the SUPG discretisation, and it can also be linked to the unphysical enstrophy dissipation at intermediate scales which affects both the Lie derivative and the flux form scheme. We also note that, when higher polynomial order is used ($h=1/86$, $r=2$), the SUPG  and Lie derivative scheme yield overshoots in the vorticity intensity. The effect is stronger for the SUPG scheme, although the overall behaviour of the solution is still well-captured.
}

\begin{figure*}[!h]
\centering
\includegraphics[scale=1.0]{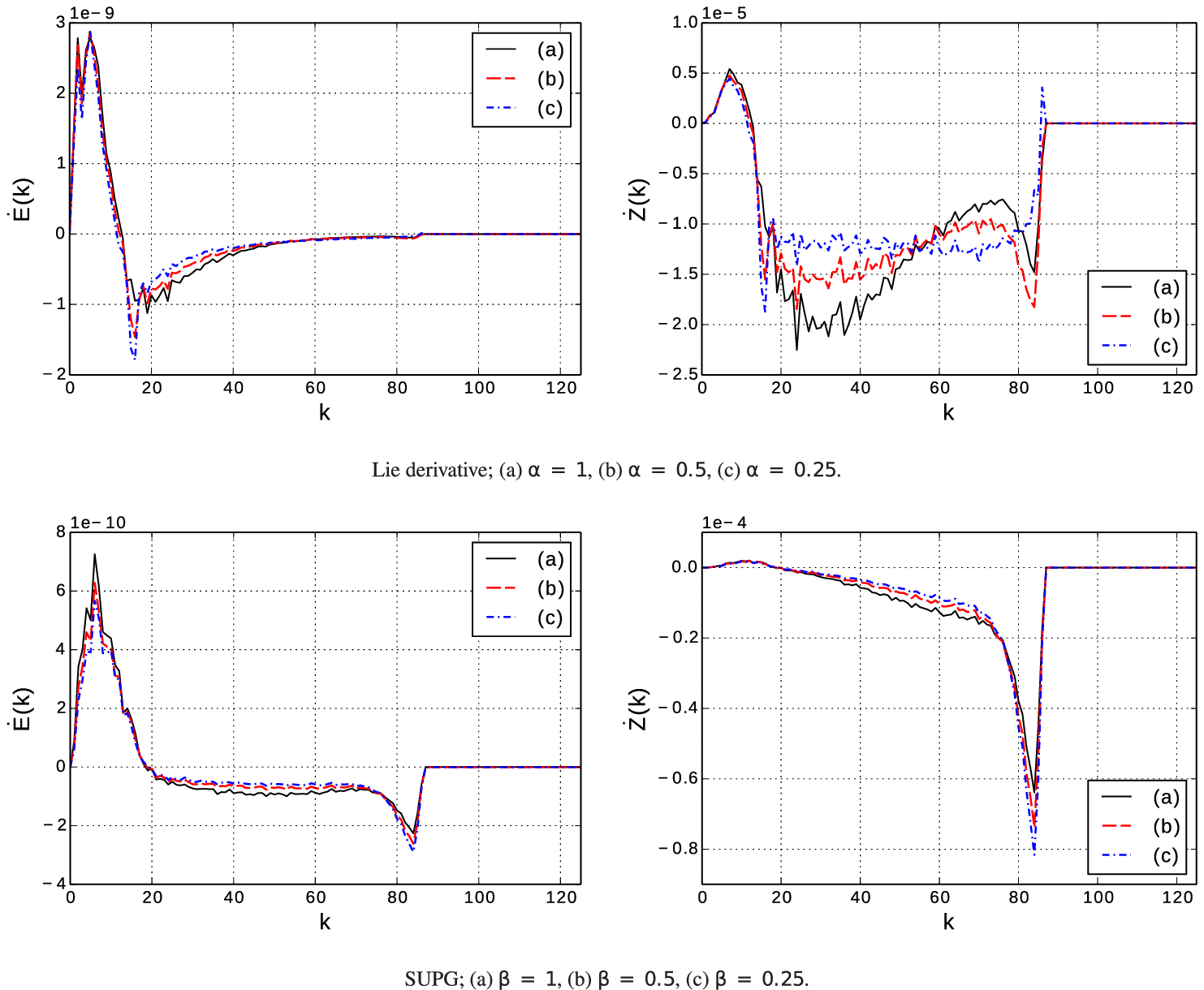}
\caption{Energy and enstrophy tendencies for the Lie derivative and SUPG scheme (with $h=1/128$, $r=1$) for different values of the stabilisation parameters $\alpha$ and $\beta$, on the spectral interval $k\leq 85$.} \label{fig:EZpar}
\end{figure*}

To complete the analysis, we examine the dependence of the results from the stabilisation parameters $\alpha$ and $\beta$, which can be used to regulate the level of dissipation. In Figure~\ref{fig:EZpar} we consider $\alpha,\beta \in\{ 0.25,0.5,1.0\}$, $h = 1/128$ and $r=1$.  Remarkably, reducing $\alpha$ in the Lie derivative scheme does not affect significantly the positive peak in the energy tendency plot. On the other hand, enstrophy dissipation moves to smaller scales. This, however, is due to the appearance of grid scale oscillations, which can be noticed in the vorticity distribution at the final time $t=350$ already for $\alpha=0.5$. Therefore, such behaviour is a signal of instability rather than improved accuracy. Similar considerations hold for the SUPG scheme. Note that for the SUPG scheme both the energy and enstrophy tendencies plots are not much affected by reducing $\beta$. However, once again, at the final time $t=350$ and for $\beta =0.25$ we observed gridscale oscillations appear in the vorticity field.

\subsection{\col{Vortex decay test}}
\begin{figure*}[!h]
\centering
\includegraphics[scale=1.0]{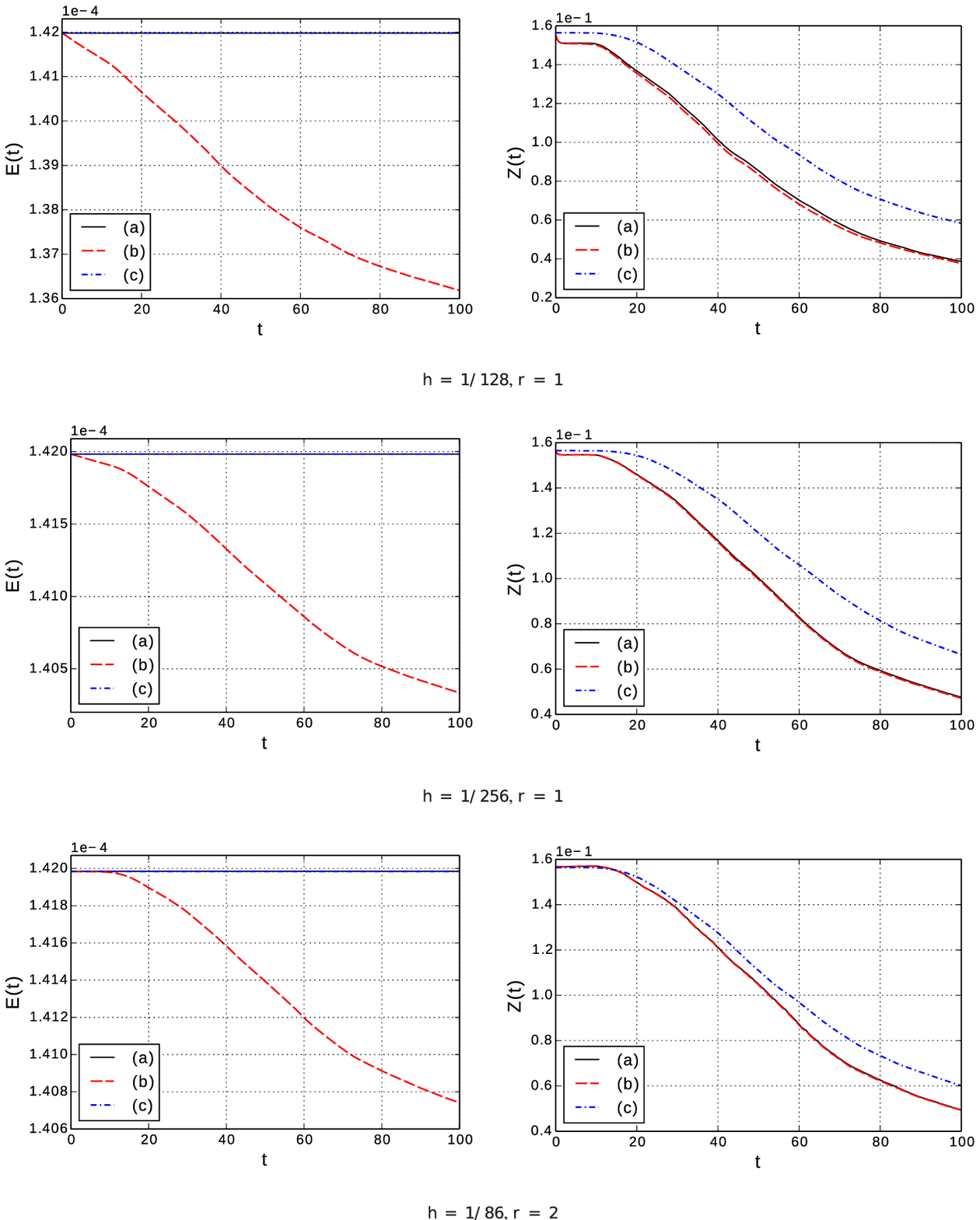}
\caption{Energy and enstrophy history for the Lie derivative (a), the flux form (b) and the SUPG (c)  scheme (with $\alpha=1$, $\beta=1$).} \label{fig:EZh}
\end{figure*}

\col{We now remove dissipation and forcing from the equations and compare the energy and enstrophy evolultion for the different schemes. We set as initial condition the one specified in Equation~\eqref{eq:ic} and evolve the solution until $t=100$, with $\alpha=1$ and $\beta=1$. The results of this test are collected in Figure~\ref{fig:EZh}. 

We observe that the Lie derivative and SUPG schemes conserve energy up to machine precision as expected. The flux form scheme dissipates energy, although this behaviour is attenuated with refinement. The enstrophy evolution for the Lie derivative and flux form scheme are very similar, and the largest differences can be seen at the lowest resolution test, i.e.\ for $h=1/128$, $r=1$, where the Lie derivative scheme is slightly less dissipative. On the other hand, the SUPG scheme dissipates substantially less enstrophy than the other two schemes. Note, in particular, the artificial drop in enstrophy at the beginning of the simulation for the Lie derivative and flux form scheme when $r=1$, which is absent from the SUPG results. 
}

\section{Conclusions}\label{sec:con}
We analysed the Lie derivative finite element discretisation of the incompressible Euler equations introduced in \cite{natale16} and the SUPG discretisation of the vorticity advection equation, in terms of energy and enstrophy tendencies in a forced turbulence test case. The main points of this paper can be summarised as follows:
\begin{itemize}
\item The Lie derivative upwind discretisation in \cite{natale16} can be reformulated to model in a simple way energy backscatter from small to large scales, by an appropriate decomposition of the velocity finite element space. The SUPG scheme, on the other hand, can be interpreted as a way to model enstrophy transfer towards unresolved scales.
\item Energy conservation in both schemes does not rely on an a priori choice of a vorticity perturbation pattern, as is usually the case for methods based on energy fixers \citep{thuburn14}.
\item The Lie derivative  scheme can be generalised to different models. This is because, on one hand, it is not exclusively a two-dimensional method; on the other, as the method descends from a Lagrangian formulation, one can derive similar algorithms by appropriately defining different energy functionals  \citep{natale16}.
\item Due to the nature of upwind stabilisation, the energy transfer is still strongly local for both schemes, however the behaviour is comparable to that obtained with standard methods combined with energy fixers \citep{thuburn14}. 
\item The SUPG scheme shows a similar qualitative behaviour to the Lie derivative discretisation but substantially better agreement with the reference solutions, probably due to higher order representation of the Jacobian.
\end{itemize}

\clearpage

\clearpage
\bibliographystyle{natbib} 
\bibliography{thesis}
 
\end{document}